\documentclass[graybox]{svmult}

\usepackage{mathptmx}       
\usepackage{helvet}         
\usepackage{courier}        
\usepackage{type1cm}        
\usepackage{makeidx}         
\usepackage{graphicx}        
\usepackage{multicol}        
\usepackage[bottom]{footmisc}

\makeindex             

\usepackage{amsmath}
\usepackage{enumitem}
\usepackage{amssymb}

\smartqed

\def\R{{\mathbb R}}  
\def\p{{\mathbb P}}  
\def\Z{{\mathbb Z}}  
\def\EE{{\mathbb E}}  

\spnewtheorem{model}[theorem]{Model}{\bf}{}
\newenvironment{Proof}[1][\unskip]{\textit{Proof #1.}}

\begin{document}
\title*{Networks, Random Graphs and Percolation}
\author{Philippe Deprez and Mario V.~W\"uthrich}
\institute{Philippe Deprez \at ETH Zurich, RiskLab,
Department of Mathematics, R\"amistrasse 101, 8092 Zurich, \email{philippe.deprez@math.ethz.ch}
\and Mario V.~W\"uthrich \at 
Swiss Finance Institute SFI Professor, 
ETH Zurich, RiskLab,
Department of Mathematics, R\"amistrasse 101, 8092 Zurich, \email{mario.wuethrich@math.ethz.ch}}
%
%
\maketitle

\abstract*{The theory of random graphs goes back to 
the late 1950s when Paul Erd\H{o}s and Alfr\'ed R\'enyi 
introduced the Erd\H{o}s-R\'enyi random graph. 
Since then many models have been developed, and 
the study of random graph models has become popular
for real-life network modelling such as social networks and
financial networks.
The aim of this overview paper is to review relevant 
random graph models for real-life network modelling. Therefore,
we analyse their properties in terms of stylised facts
of real-life networks.}

\abstract{The theory of random graphs goes back to 
the late 1950s when Paul Erd\H{o}s and Alfr\'ed R\'enyi 
introduced the Erd\H{o}s-R\'enyi random graph. 
Since then many models have been developed, and 
the study of random graph models has become popular
for real-life network modelling such as social networks and
financial networks.
The aim of this overview is to review relevant 
random graph models for real-life network modelling. Therefore,
we analyse their properties in terms of stylised facts
of real-life networks.}

\section{Stylised facts\index{stylised facts} of real-life networks}
\label{section introduction}
A network is a set of particles that may be linked to each other. The particles represent
individual network participants and the links illustrate how they interact among each other, for 
an example see Fig.~\ref{ER_and_NSW} below. 
Such networks appear in many real-life situations, for instance, there are virtual social 
networks with different users that
communicate with (are linked to) each other, see Newman et al.~\cite{NSW2}, 
or there are financial networks such as the banking system that exchange lines of credits
with each other,
see Amini et al.~\cite{Cont1} and Cont et al.~\cite{Cont2}. These two examples 
represent rather recently established real-life networks 
that originate from new technologies and industries but, of course, the study of network
models is much older motivated by studies in sociology or questions about interacting
particle systems in physics.

Such real-life networks, in particular
social networks, have been studied on many different empirical data sets. These studies have
raised several stylised facts about large real-life networks
 that we would briefly like to enumerate, for more details
see Newman et al.~\cite{NSW2} and Section 1.3 in Durrett \cite{Durrett} and the references therein.
\begin{enumerate}
\item Many pairs of distant particles are connected by a very short chain of links. This is
sometimes called the ``small-world'' effect\index{small-world effect}. 
Another interpretation of the small-world effect is
the observation that the typical distance of any two particles in real-life networks is at most six links,
see Watts \cite{Watts} and Section 1.3 in Durrett \cite{Durrett}. 
The work of Watts \cite{Watts} was inspired by the statement of his father saying that ``he 
is only six handshakes away from the president of the United States''.
For other interpretations of the
small-world effect we refer to Newman  et al.~\cite{NSW2}.
\item The clustering property\index{clustering property} of real-life networks is often observed which means that linked particles
tend to have common friends. 
\item The distribution of the number of links\index{degree!distribution} of a single particle is heavy-tailed, i.e.~its
survival probability has a power law decay. In many real-life networks 
the power law constant (tail parameter) $\tau$
is estimated between 1 and 2 (finite mean and infinite variance, see also (\ref{NSW degree}) below). Section
1.4 in Durrett \cite{Durrett} presents the following examples:
\begin{itemize}
\item number of oriented links on web pages: $\tau \approx 1.5$,
\item routers for e-mails and files: $\tau \approx 1.2$,
\item movie actor network: $\tau \approx 1.3$,
\item citation network Physical Review D: $\tau \approx 1.9$.
\end{itemize}
Typical real-life networks are heavy-tailed in particular if maintaining links is free of costs.
\end{enumerate}
Since real-life networks are too complex to be modelled particle by particle and link by link, 
researchers have
developed many models in random graph theory that help to understand the geometry of 
such real-life networks. The aim of this overview paper is to review relevant models
in random graph theory, in particular,
we would like to analyse whether these models fulfil the stylised facts.
Standard literature on random graph and percolation theory is Bollob\'{a}s \cite{Bollobas},
Durrett \cite{Durrett}, Franceschetti-Meester \cite{Franceschetti}, 
Grimmett \cite{Grimmett1, Grimmett2} and Meester-Roy \cite{Meester}.

\section{Erd\H{o}s-R\'enyi random graph\index{random graph!Erd\H{o}s-R\'enyi}}
\begin{figure}[htb!]
\begin{center}
\begin{minipage}[t]{0.45\textwidth}
\begin{center}
\includegraphics[width=.9\textwidth]{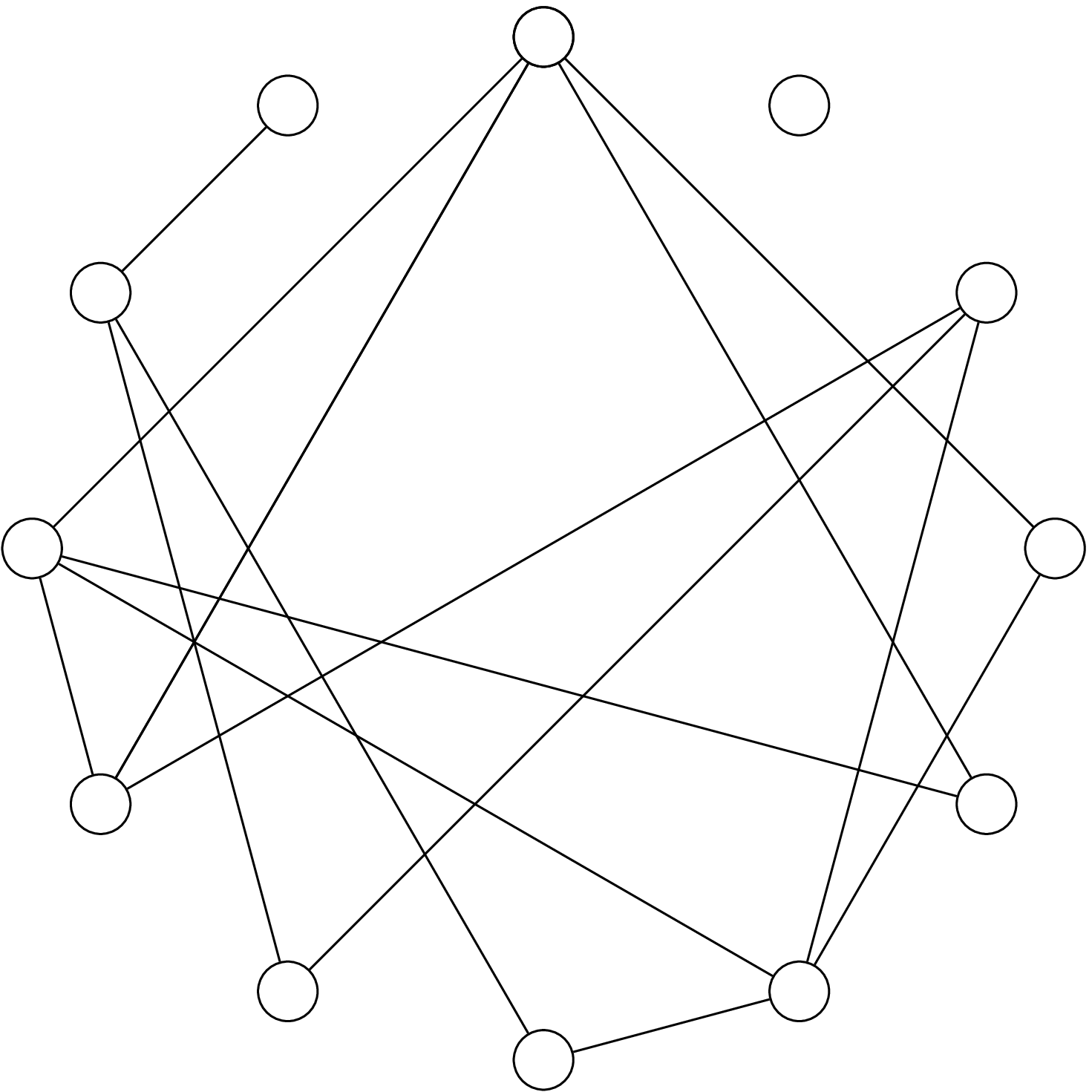}
\end{center}
\end{minipage}
\begin{minipage}[t]{0.45\textwidth}
\begin{center}
\includegraphics[width=.9\textwidth]{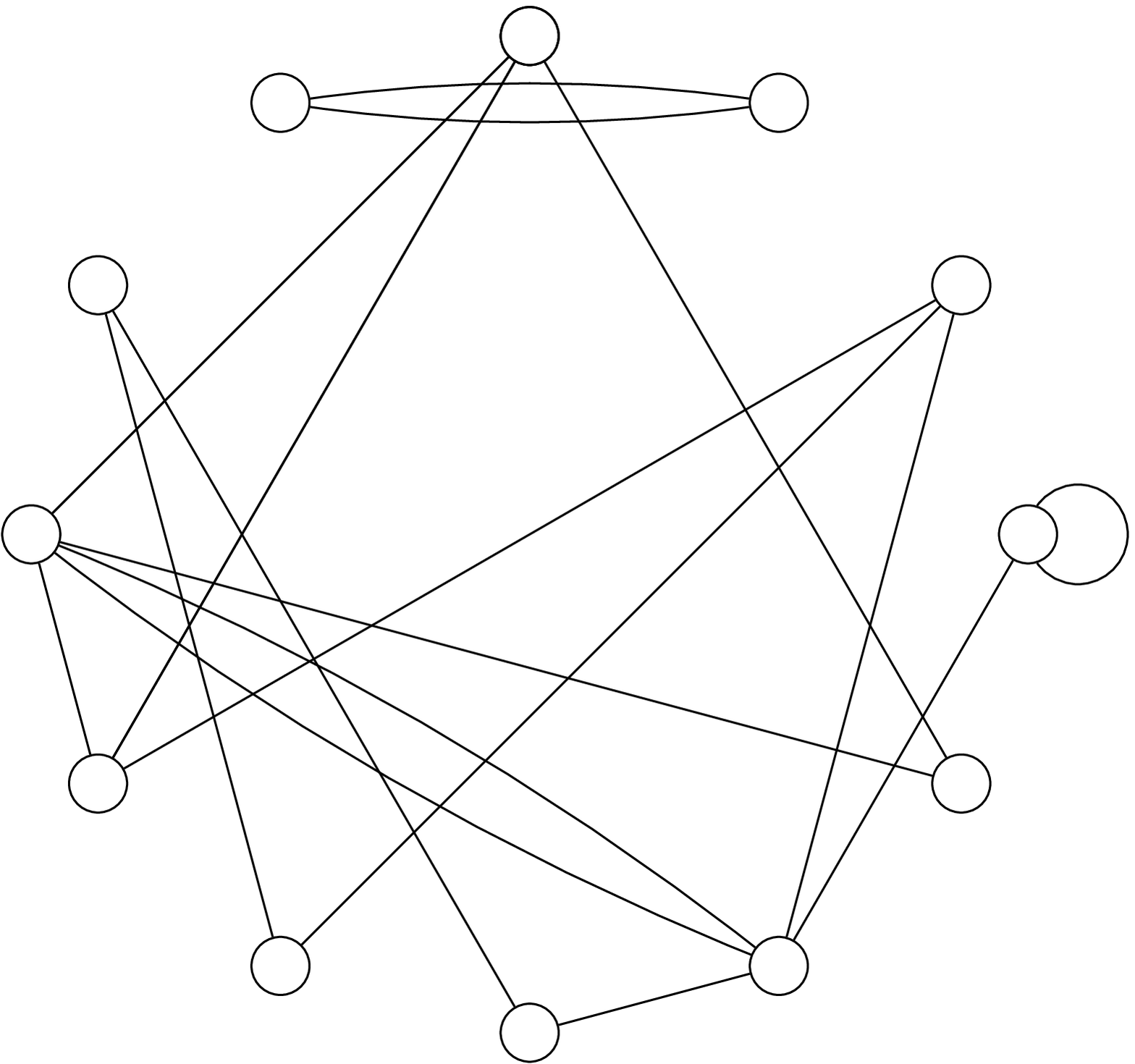}
\end{center}
\end{minipage}
\end{center}
\caption{{\bf lhs}: ER random graph; {\bf rhs}: NSW random graph} \label{ER_and_NSW}
\end{figure}

We choose a set of particles $V_n=\{1,\ldots, n\}$ 
for fixed $n\in \mathbb{N}$. Thus, $V_n$ contains $n$ particles.
The Erd\H{o}s-R\'enyi (ER) random graph introduced in the late 1950s, see \cite{ER},
attaches to every pair of particles $x,y\in V_n$, $x\neq y$, 
independently an edge with fixed probability\index{edge probability} $p \in (0,1)$, i.e.,
\begin{equation}\label{Erdos-Renyi}
\eta_{x,y}=\eta_{y,x} =
\left\{ 
\begin{array}{ll}
1 \qquad & \text{ with probability $p$,}\\
0 \qquad & \text{ with probability $1-p$,}
\end{array}
\right.
\end{equation}
where $\eta_{x,y}=1$ means that there is an edge between $x$ 
and $y$, and $\eta_{x,y}=0$ means that there is {\it no} edge between $x$ 
and $y$. Identity $\eta_{x,y}=\eta_{y,x}$ illustrates that
we have an undirected random graph.
We denote this random graph model by ${\rm ER}(n,p)$.
In Fig.~\ref{ER_and_NSW} (lhs) we provide an example for $n=12$, observe that
this realisation of the
ER random graph has one isolated particle and the remaining ones lie in the
same connected component.

We say that $x$ and $y$ are {\it adjacent}\index{adjacent}
if $\eta_{x,y}=1$. We say that $x$ and $y$ are {\it connected}\index{connected}
if there exists a path of adjacent particles from $x$ to $y$. 
We define 
the {\it degree}\index{degree} ${\cal D}(x)$ of particle $x$ to 
be the number of adjacent particles of $x$ in $V_n$. 
Among others, general random graph theory is concerned with the limiting behaviour of
the ER random graph ${\rm ER}(n,p_n)$ for $p_n=\vartheta/n$, $\vartheta>0$,
as $n\to\infty$. Observe that
for $k\in \{0,\ldots, n-1\}$ we have, see for instance Lemma 2.9 in \cite{W},
\begin{equation}\label{ER degree}
g_k=g^{(n)}_k=\p\left[{\cal D}(x)=k\right]
=\binom{n-1}{k} p_n^k \left(1-p_n\right)^{n-1-k} \quad \rightarrow \quad
e^{-\vartheta} \frac{\vartheta^k}{k!},
\end{equation}
as $n\to \infty$.
We see that the degree distribution  
of a fixed particle $x \in V_n $ with edge
probability $p_n$ converges for $n\to \infty$ to a Poisson
distribution with parameter $\vartheta>0$. In particular, this limiting distribution
is light-tailed and, therefore, the ER graph does not fulfil the stylised fact
of having a power law decay of the degree distribution.

The ER random graph has a phase transition\index{phase transition} at $\vartheta=1$, reflecting different regimes for the
size of the largest connected component in the ER random graph.
For $\vartheta<1$, all connected components\index{connected component} are small, the largest
being of order $\mathcal{O}(\log n)$, as $n\to \infty$. 
For $\vartheta >1$, there is a constant $\chi(\vartheta)>0$ and the largest connected component
of the ER random graph is
of order $\chi(\vartheta)n$, as $n\to\infty$, 
and all other connected components are small,
see Bollob\'{a}s \cite{Bollobas} and 
Chapter 2 in Durrett \cite{Durrett}.
At criticality ($\vartheta=1$) the largest connected component 
is of order $n^{2/3}$, however, this
analysis is rather sophisticated, see Section 2.7 in Durrett \cite{Durrett}.

Moreover, the ER random graph has only very few
complex connected components such as cycles (see Section 2.6 in Durrett \cite{Durrett}):
for $\vartheta \neq 1$ most connected components are trees, only a few 
connected components have triangles and cycles, and
only the largest connected component (for $\vartheta>1$) is more complicated. 
At criticality the situation is more complex, a few large connected components emerge and finally
merge to the largest connected component as $n\to \infty$.

\section{Newman-Strogatz-Watts random graph\index{random graph!Newman-Strogatz-Watts}}
\label{Section NSW graph}
The approach of Newman-Strogatz-Watts (NSW) \cite{NSW1, NSW2} aims
at directly describing the degree distribution $(g_k)_{k \ge 0}$ of 
${\cal D}(x)$ for a given particle $x \in V_n$
($n\in \mathbb{N}$ being large).
The aim is to modify the degree distribution in
(\ref{ER degree}) so that we obtain a power law
distribution.
Assume that any particle $x \in V_n$ has a
degree distribution\index{degree!distribution} of the form $g_0=0$ and 
\begin{equation}\label{NSW degree}
g_k=\p\left[{\cal D}(x)=k\right] \sim c k^{-(\tau+1)},
\qquad \text{ as $k\to \infty$,}
\end{equation}
for given tail parameter $\tau>0$ and $c>0$.
Note that $\sum_{k\ge 1}k^{-(\tau+1)}<1+1/\tau$ which implies that $c>0$ is admissible.
By definition the survival probability of 
this degree distribution has a power law with tail parameter $\tau>0$. However,
this choice (\ref{NSW degree})
does not explain how one obtains an explicit graph from the degrees ${\cal D}(x)$,
$x\in V_n$. The graph construction is done 
by the Molloy-Reed \cite{MR} algorithm\index{Molloy-Reed algorithm}: 
attach
to each particle $x \in V_n$ exactly ${\cal D}(x)$ ends of edges and then choose these ends randomly
in pairs (with a small modification if the total number of ends is odd). 
This will provide a random graph with the desired degree distribution. 
In Fig.~\ref{ER_and_NSW} (rhs) we provide an example for $n=12$, observe that
this realisation of the
NSW random graph has two connected components. The Molloy-Reed construction may provide multiple edges
and self-loops, but if ${\cal D}(x)$ has finite second moment ($\tau>2$) then there are only a few multiple
edges and self-loops, as $n\to \infty$, see
Theorem 3.1.2 in Durrett \cite{Durrett}. However, in view of real-life networks we are rather
interested into tail parameters $\tau \in (1,2)$ for which we so far have no control on
multiple edges and self-loops.

Newman et al.~\cite{NSW1, NSW2} have analysed this random graph by basically considering cluster growth
in a two-step branching process. Define the probability generating function 
of the first generation by
\begin{equation*}
G_0(z) = \EE \left[z^{{\cal D}(x)}\right] = \sum_{k \ge 1} g_k z^k,\qquad
\text{ for $z\in \R$.}
\end{equation*}
Note that we have $G_0(1)=1$ and $\mu=\EE[{\cal D}(x)]=G'_0(1)$ (supposed that the latter exists). 
The second generation has then probability generating function given by
\begin{equation*}
G_1(z) = \sum_{k\ge 0} \frac{(k+1)g_{k+1}}{\mu}~ z^{k}
= \sum_{k\ge 1} \frac{kg_k}{\mu}~ z^{k-1},\qquad
\text{ for $z\in \R$,}
\end{equation*}
where the probability weights are specified by $\widetilde{g}_k=(k+1)g_{k+1}/\mu$ for $k\ge 0$.
For $\tau>2$ the second generation has finite mean given by
\begin{equation*}
\vartheta= \sum_{k\ge 0} k \widetilde{g}_k
= \sum_{k\ge 0} k~\frac{(k+1)g_{k+1}}{\mu}
=\frac{1}{\mu} \sum_{k\ge 1} (k-1)kg_k.
\end{equation*}
Note that the probability generating functions are related to each other by
$G_0'(z)= \mu G_1(z) = G_0'(1)G_1(z)$.
Similar to the ER random graph there is a phase transition\index{phase transition}
in this model. It is determined by
the mean $\vartheta$ of the second generation, see (5)--(6) in Newman et al.~\cite{NSW2} 
and Theorems 3.1.3 and 3.2.2 in Durrett \cite{Durrett}: 
for $\vartheta>1$ the largest connected component has size of order $\chi(\vartheta)n$, 
as $n\to \infty$. The fraction $\chi(\vartheta)=1-G_0(z_0)$ 
is found by choosing $z_0$ to be the smallest fixed point of $G_1$ in $[0,1]$. Moreover, no
other connected component has size of order larger than $\mathcal{O}(\log n)$.
Note that we require finite variance $\tau>2$ for $\vartheta$ to exist.

If $\vartheta <1$ the distribution of the size of the connected component of a fixed
particle converges in distribution to a limit with mean $1+\mu/(1-\vartheta)$, as $n\to \infty$,
see Theorem 3.2.1 in Durrett \cite{Durrett}. The size of the largest connected component
in this case ($\tau>2$ and $\vartheta<1$) is conjectured to be of order $n^{1/\tau}$:
the survival probability of the degree distribution has asymptotic behaviour 
of order $k^{-\tau}$, therefore the largest degree of $n$ independent degrees
has size of order $n^{1/\tau}$, which leads to the same conjecture for the largest
connected component, see also Conjecture 3.3.1 in Durrett \cite{Durrett}.

From a practical point of view the interesting regime is $1<\tau<2$ because many 
real-life networks have such a tail behaviour, see Section
1.4 in Durrett \cite{Durrett}. In this case we have $\vartheta=\infty$
and an easy consequence is that the largest connected component grows proportionally
to $n$ (because this model dominates a model with finite second moment and mean
of the second generation being bigger than 1). In this regime $1<\tau<2$ we can
study the graph distance\index{graph distance} of two randomly chosen particles 
(counting the number of edges connecting them) 
in the largest connected component,
see Section 4.5 in Durrett \cite{Durrett}. In the Chung-Lu model \cite {CL0, CL},
which uses a variant to the Molley-Reed \cite{MR} algorithm, it is proved 
that this graph distance behaves as ${\mathcal O}(\log\log n)$, see Theorem 4.5.2
in Durrett \cite{Durrett}. Van der Hofstadt et al.~\cite{Remco2} 
obtain the same asymptotic behaviour ${\mathcal O}(\log\log n)$
for the NSW random graph in the case 
$1<\tau<2$. Moreover, in their Theorem 1.2 \cite{Remco2} they also state that
this graph distance behaves as ${\mathcal O}(\log n)$ for $\tau>2$.
These results on the graph distances can be interpreted as the small-world
effect because two randomly chosen particles in $V_n$ are connected by
very few edges.

We conclude that  NSW random graphs have heavy tails for the degree distribution
choices according to (\ref{NSW degree}). Moreover, the graph distances
have a behaviour that can be interpreted as small-world effect.

Less desirable features of NSW random graphs are that they may have
self-loops and multiple edges. Moreover, the NSW random graph is expected to be
locally rather sparse leading to locally tree-like structures, see also 
Hurd-Gleeson \cite{Hurd1}. That is, we do 
not expect to get a reasonable local graph geometry and the required clustering property.
Variations considered allowing for statistical interpretations in terms
of likelihoods include the works of Chung-Lu \cite{CL0, CL} and Olhede-Wolfe \cite{OW}.

\section{Nearest-neighbour bond percolation\index{percolation!nearest-neighbour bond}}
\label{nearest-neighbour bond percolation}
In a next step we would like to embed the previously introduced random graphs and the corresponding particles into
Euclidean space. This will have the advantage of obtaining 
a natural distance function between particles, and it will allow to compare
Euclidean distance to graph distance between particles (counting 
the number
of edges connecting two distinct particles). 
Before giving the general random graph model we restrict
ourselves to the nearest-neighbour bond percolation model on the lattice $\Z^d$ because
this model is the basis for many derivations. 
More general and flexible 
random graph models are provided in the subsequent sections.

Percolation theory was first presented by Broadbent-Hammersley \cite{BH}. 
It was mainly motivated
by questions from physics, but these days percolation models are recognised to be very useful in several
fields. Key monographs on nearest-neighbour
bond percolation theory are Kesten \cite{Kesten} and Grimmett \cite{Grimmett1, Grimmett2}.

\begin{figure}[htb!]
\begin{center}
\begin{minipage}[t]{0.45\textwidth}
\begin{center}
\includegraphics[width=.9\textwidth]{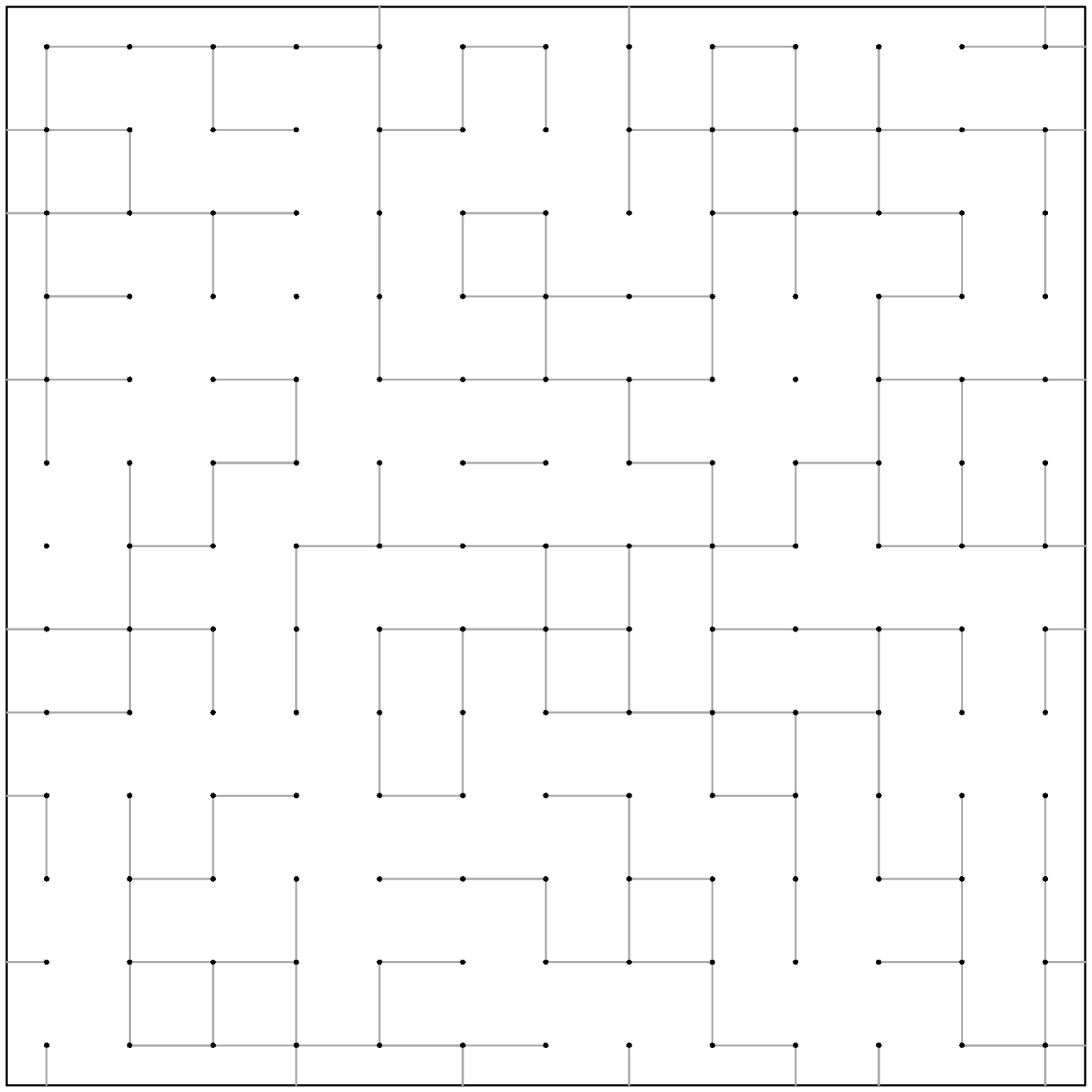}
\end{center}
\end{minipage}
\begin{minipage}[t]{0.45\textwidth}
\begin{center}
\includegraphics[width=.9\textwidth]{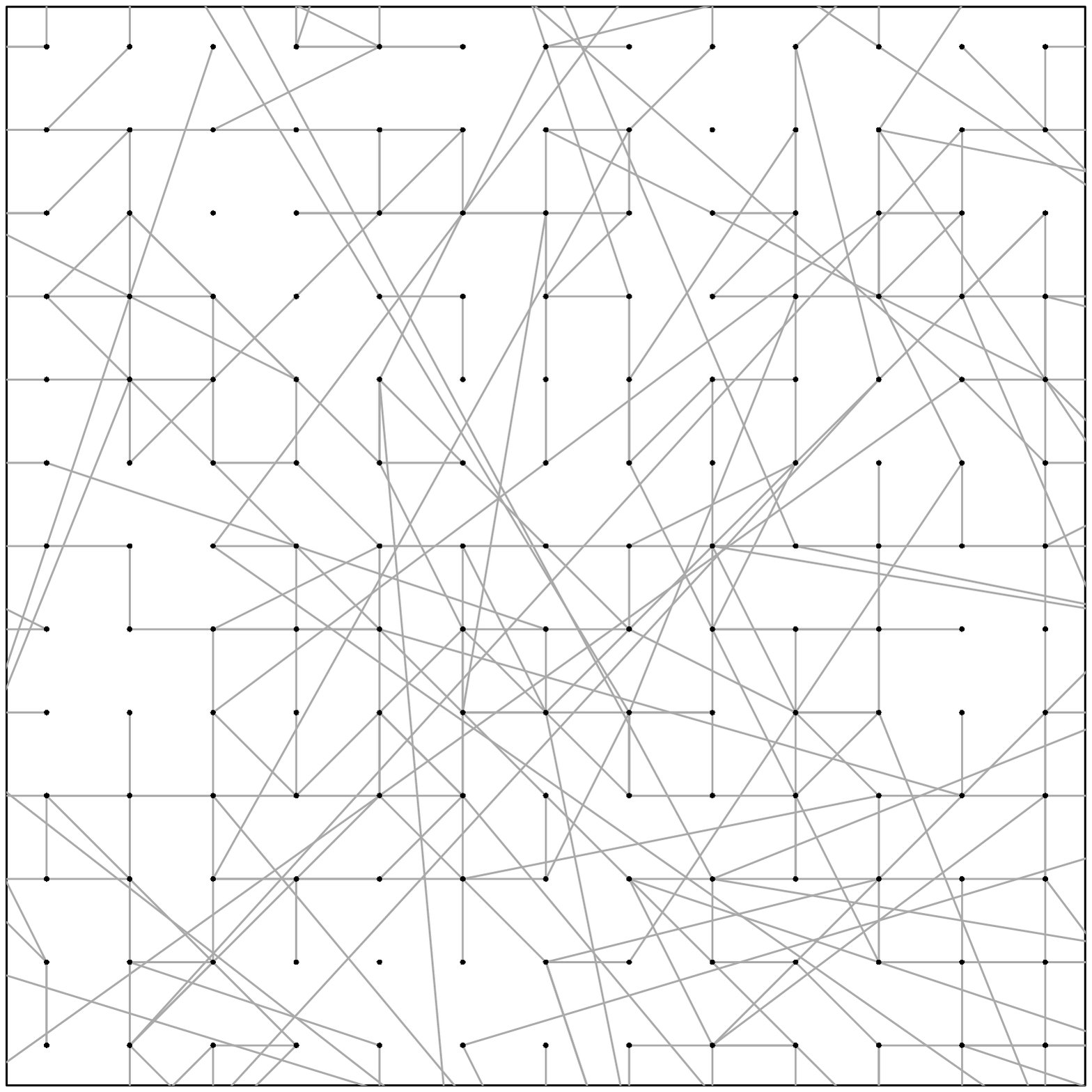}
\end{center}
\end{minipage}
\end{center}
\caption{{\bf lhs}: nearest-neighbour percolation; 
{\bf rhs}: homogeneous long-range percolation} \label{NNBP_and_HomLRP}
\end{figure}

Choose a fixed dimension $d\in \mathbb{N}$ and consider the square lattice $\Z^d$. 
The vertices of this square lattice are the particles and we say that two particles $x,y\in \Z^d$ are 
nearest-neighbour particles if $\|x-y\|=1$ (where $\|\cdot\|$ denotes the Euclidean norm).
We attach at random edges to nearest-neighbour particles $x,y \in \Z^d$, independently of all
other edges, with a fixed edge probability\index{edge probability} $p\in [0,1]$, that is,
\begin{equation}\label{percolation}
\eta_{x,y}=\eta_{y,x} =
\left\{ 
\begin{array}{ll}
1_{\{\|x-y\|=1\}} \qquad & \text{ with probability $p$,}\\
0 \qquad & \text{ with probability $1-p$,}
\end{array}
\right.
\end{equation}
where $\eta_{x,y}=1$ means that there is an edge between $x$ 
and $y$, and $\eta_{x,y}=0$ means that there is {\it no} edge between $x$ 
and $y$. The resulting graph is called nearest-neighbour (bond) 
random graph  in $\Z^d$, 
see Fig.~\ref{NNBP_and_HomLRP} (lhs) for an illustration. Two
particles $x, y \in \Z^d$ are connected if there exists a path of nearest-neighbour edges 
connecting $x$ and $y$.
It is immediately clear that this random graph does not fulfil the small-world 
effect because one
needs at least $\|x-y\|$ edges to connect $x$ and $y$, i.e.~the number of edges grows at
least linearly in the Euclidean distance between particles $x,y\in \Z^d$. The degree 
distribution 
is finite because there are at most $2^d$ nearest-neighbour edges, more precisely, the degree has a binomial
distribution with parameters $2^d$ and $p$. 
 We present this square lattice model because it is an interesting basis for the
development of more complex models. Moreover, this model is at the heart of many proofs
in percolation problems which are based on so-called renormalisation techniques, 
see Sect.~\ref{Renormalisation techniques} below for a concrete example. 

In percolation theory, the object of main
interest is the connected component\index{connected component} of a given particle $x\in \Z^d$
which we denote by
\begin{equation*}
{\cal C}(x) = \left\{ y \in \Z^d:~\text{$x$ and $y$ are connected by a path of nearest-neighbour edges}\right\}.
\end{equation*}
By translation invariance it suffices to define the percolation 
probability\index{percolation!probability} at the origin
\begin{equation*}
\theta(p)= \p_p \left[ |{\cal C}(0)|=\infty \right],
\end{equation*}
where $|{\cal C}(0)|$ denotes the size of the connected component of the origin and
$\p_p$ is the product measure on the  possible nearest-neighbour edges with edge probability $p\in [0,1]$,
see Grimmett \cite{Grimmett1}, Section 2.2. 
The {\it critical probability}\index{critical probability} $p_c=p_c(\Z^d)$
is then defined by
\begin{equation*}
p_c = \inf \left\{ p\in (0,1]:~\theta(p)>0\right\}.
\end{equation*}
Since the percolation probability $\theta(p)$ is non-decreasing, the critical probability
is well-defined. We have the following result, see Theorem 3.2 in Grimmett \cite{Grimmett1}.
\begin{theorem} \label{bernoulli theo 1}
For nearest-neighbour bond percolation in $\Z^d$ we have
\begin{itemize}
\item[(a)] for $d=1$: $p_c(\Z)=1$; and
\item[(b)] for $d\ge 2$: $p_c(\Z^d) \in (0,1)$.
\end{itemize}
\end{theorem}
This theorem says that there is a non-trivial phase 
transition\index{phase transition} in $\Z^d$, $d\ge 2$. This needs to be
considered together with the following result which goes back to Aizenman et al.~\cite{AKN},
Gandolfi et al.~\cite{GGR} and Burton-Keane \cite{BK}. Denote by ${\cal I}$ the number of infinite 
connected components. Then we have the following statement, see Theorem 7.1 in Grimmett \cite{Grimmett1}.
\begin{theorem}\label{bernoulli theo 2}
For any $p\in (0,1)$ either $\p_p[{\cal I}=0]=1$ or $\p_p[{\cal I}=1]=1$.
\end{theorem}
Theorems \ref{bernoulli theo 1} and \ref{bernoulli theo 2} imply that there is a {\it unique} infinite
connected component for $p>p_c(\Z^d)$, a.s.
This motivates the notation ${\cal C}_\infty$ for the unique infinite connected 
component\index{connected component!infinite}
for the given edge configuration $(\eta_{x,y})_{x,y}$
in the case  $p>p_c(\Z^d)$. ${\cal C}_\infty$ may be considered as an infinite (nearest-neighbour) network 
on the particle system $\Z^d$
and we can study its geometrical and topological properties. 
Using a duality argument, Kesten \cite{Kesten12} proved that $p_c(\Z^2)=1/2$ and monotonicity
then provides $p_c(\Z^{d+1})\le p_c(\Z^{d})\le p_c(\Z^{2})=1/2$ for $d\ge 2$.

One object
of interest is the so-called {\it graph distance} 
(chemical distance)\index{graph distance}\index{chemical distance} between $x,y\in \Z^d$,
which is for a given edge configuration  defined by
\begin{eqnarray*}
d(x,y)&=& \text{minimal length of path connecting $x$ and $y$ by}
\\&&\hspace{3cm} \text{nearest-neighbour edges 
$\eta_{z_1,z_2}=1$},
\end{eqnarray*}
where this is defined to be infinite if there is no nearest-neighbour path connecting
$x$ and $y$ for the given edge configuration. We have already mentioned that
$d(x,y)\ge \|x-y\|$ because this is the minimal number of nearest-neighbour edges
we need to cross from $x$ to $y$.
Antal-Pisztora \cite{AP} have proved the following upper bound.
\begin{theorem} \label{Antal theorem}
Choose $p>p_c(\Z^d)$. There exists a positive constant $c=c(p,d)$ such
that, a.s.,
\begin{equation*}
\limsup_{\|x\|\to\infty} \frac{1}{\|x\|}d(0,x)1_{\{\text{\rm 0 and $x$ are connected}\}} \le c.
\end{equation*}
\end{theorem}

\section{Homogeneous long-range percolation}\index{percolation!homogeneous long-range}
\label{homogeneous long-range percolation}
Long-range percolation is the
first extension of  nearest-neighbour bond percolation. It allows for 
edges between any pair of particles $x,y\in \Z^d$. Long-range percolation was originally
introduced by Schulman \cite{Schulman} in one dimension. Existence and
uniqueness of the infinite connected component in long-range
percolation was proved by Schulman \cite{Schulman} and
Newman-Schulman \cite{NS} for $d=1$ and by Gandolfi et al.~\cite{GKN} for $d\ge 2$.

Consider again the percolation model on the lattice $\Z^d$, but we now choose the edges
differently. Choose $p\in [0,1]$, $\lambda>0$ and $\alpha>0$  fixed and define the edge
probabilities\index{edge probability} for $x,y\in \Z^d$ by
\begin{equation}\label{long-range percolation}
p_{x,y} =
\left\{ 
\begin{array}{ll}
p \qquad & \text{ if $\|x-y\|=1$,}\\
1-\exp(-\lambda \|x-y\|^{-\alpha}) \qquad & 
\text{ if $\|x-y\|>1$.}
\end{array}
\right.
\end{equation}
Between any pair $x,y\in \Z^d$ we attach an edge, independently of all other
edges, as follows
\begin{equation*}
\eta_{x,y}=\eta_{y,x} =
\left\{ 
\begin{array}{ll}
1 \qquad & \text{ with probability $p_{x,y}$,}\\
0 \qquad & \text{ with probability $1-p_{x,y}$.}
\end{array}
\right.
\end{equation*}
We denote the resulting product measure on the edge configurations by $\p_{p,\lambda, \alpha}$. 
Figure \ref{NNBP_and_HomLRP} (rhs) shows part of a realised configuration. 
We say that the particles $x$ and $y$ are {\it adjacent} if there is an edge $\eta_{x,y}=1$
between $x$ and $y$. We say that $x$ and $y$ are {\it connected}
if there exists a path of adjacent particles in $\Z^d$ that connects $x$ and $y$. The connected component\index{connected component}
of $x$ is given by 
\begin{equation*}
{\cal C}(x)=\left\{ y \in \Z^d:~\text{$x$ and $y$ are connected} \right\}.
\end{equation*}
We remark that the edge probabilities $p_{x,y}$ used in the literature have a more general form.
Since for many results only the asymptotic behaviour of $p_{x,y}$ as $\|x-y\|\to \infty$ is 
relevant, we have decided to choose the explicit (simpler) 
form (\ref{long-range percolation}) because
this also fits to our next models. Asymptotically we have  the following power law
\begin{equation*}
p_{x,y}~ \sim~ \lambda \|x-y\|^{-\alpha}, \qquad \text{ as $\|x-y\|\to \infty$.}
\end{equation*}
Theorem \ref{bernoulli theo 1} (b) immediately implies that we have percolation in $\Z^d$,
$d\ge 2$, for $p$ sufficiently close to 1. We have the following theorem, see
Theorem 1.2 in Berger \cite{Berger1}.

\begin{theorem} \label{theorem infinite cluster}
For long-range percolation in $\Z^d$ we have, in an a.s.~sense,\index{phase transition}
\begin{itemize}
\item[(a)] for $\alpha \le d$: there is an infinite connected component;
\item[(b)] for $d\ge 2$ and $\alpha > d$: for $p$ sufficiently close to 1 there
is an infinite connected component;
\item[(c)] for $d=1$:
\begin{itemize}
\item[(1)] $\alpha>2$: there is no infinite connected component;
\item[(2)] $1<\alpha<2$: for $p$ sufficiently close to 1 there
is an infinite connected component;
\item[(3)] $\alpha=2$ and $\lambda>1$: for $p$ sufficiently close to 1 there
is an infinite connected component;
\item[(4)] $\alpha=2$ and $\lambda\le 1$: there is no infinite connected component.
\end{itemize}
\end{itemize} 
\end{theorem}
The case $\alpha\le d$ follows from an infinite degree distribution for a given particle, 
i.e.~for $\alpha \le d$ we have, a.s.,
\begin{equation}\label{homogeneous lattice}
{\cal D}(0)=\left|\left\{ x \in\Z^d:~\text{$0$ and $x$ are adjacent}\right\}\right|=\infty,
\end{equation}
and for $\alpha>d$ the degree distribution is light-tailed (we give a proof in the
continuum space model in Sect.~\ref{homogeneous Poisson model}, because the proof turns
out to be straightforward in continuum space).
Interestingly, we now also obtain a non-trivial phase transition in the one dimensional
case $d=1$ once
long-range edges are sufficiently likely, i.e.~$\alpha$ is sufficiently small. At 
criticality $\alpha=2$ also the decay scaling constant $\lambda>0$ matters. The case $d\ge 2$
is less interesting because it is in line with nearest-neighbour bond percolation. The main
interest of adding long-range edges is the study of the resulting geometric properties of connected 
components ${\cal C}(x)$.
We will state below that there are three different regimes:
\begin{itemize}
\item $\alpha \le d$ results in an infinite degree distribution, a.s., see (\ref{homogeneous lattice});
\item $d < \alpha < 2d$ has finite degrees but  is still in the regime of small-world behaviour;
\item $\alpha >2d$ behaves as nearest-neighbour bond percolation.
\end{itemize}
We again focus on the graph distance 
\begin{equation}\label{definition of graph distance d}
d(x,y)= \text{minimal number of edges that connect $x$ and $y$},
\end{equation}
where this is defined to be infinite if $x$ and $y$ do not belong to the same connected
component, i.e.~$y \notin {\cal C}(x)$.
For $\alpha<d$ we have infinite degrees and the infinite connected component ${\cal C}_\infty$
contains
all particles of $\Z^d$, a.s. Moreover, Benjamini et al.~\cite{BKPS} prove in Example 6.1 that
the graph distance is bounded, a.s., by
\begin{equation*}
\left\lceil \frac{d}{d-\alpha} \right\rceil .
\end{equation*}
The case $\alpha \in (d,2d)$ is considered in Biskup \cite{Biskup}, Theorem 1.1, and
in Trapman \cite{Trapman}. They have proved the following result: 
\begin{theorem}\label{theorem homogeneous 1}
Choose $\alpha \in (d,2d)$ and assume, a.s., that there exists a unique infinite connected
component ${\cal C}_\infty$. Then for all $\epsilon>0$ we have
\begin{equation*}
\lim_{\|x\|\to\infty} \p_{p,\lambda,\alpha}
\left[\left.\Delta -\epsilon \le \frac{\log d(0,x)}{\log \log \|x\|} \le
\Delta +\epsilon\right|0,x\in {\cal C}_\infty\right]=1,
\end{equation*}
where $\Delta^{-1}=\log_2 (2d/\alpha)$.
\end{theorem}
This result says that the graph distance $d(0,x)$ is roughly of order
$(\log \|x\|)^\Delta$ with $\Delta=\Delta(\alpha, d)>1$. Unfortunately, the known bounds
are not sufficiently sharp to give more precise asymptotic statements.
Theorem \ref{theorem homogeneous 1} can be interpreted as small-world effect since
it tells us that long Euclidean distances can be crossed by a few edges. 
For instance, $d=2$ and
$\alpha=2.5$ provide $\Delta= 1.47$ and we get $(\log \|x\|)^\Delta=26.43$ for $\|x\|=10,000$, i.e.~a Euclidean distance of 10,000 is crossed in roughly 26 edges.

The case $\alpha>2d$ is considered in Berger \cite{Berger2}.
\begin{theorem} \label{theorem homogeneous 2}
If $\alpha>2d$ we have, a.s.,
\begin{equation*}
\liminf_{\|x\|\to \infty} \frac{d(0,x)}{\|x\|}>0.
\end{equation*}
\end{theorem}
This result proves that for $\alpha>2d$ the graph distance behaves as
in nearest-neighbour bond percolation, because it grows linearly in $\|x\|$.
The proof of an upper bound is still open, but we expect a result similar
to Theorem \ref{Antal theorem} in nearest-neighbour bond percolation, 
see Conjecture 1 of Berger \cite{Berger2}. 

We conclude that this model has a small-world effect for $\alpha <2d$. It also has some kind of clustering
property because particles that are close share an edge more commonly, which gives a structure that is
locally more dense, see Corollary 3.4 in Biskup
\cite{Biskup}. But the degree distribution is light-tailed which motivates to extend
the model by an additional ingredient. This is done in the next section.

\section{Heterogeneous long-range percolation\index{percolation!heterogeneous long-range}}
Heterogeneous long-range percolation extends the previously introduced long-range percolation
models on the lattice $\Z^d$. Deijfen et al.~\cite{Remco} have introduced this model  under the name of scale-free percolation\index{percolation!scale-free}. The idea is to place additional weights\index{weights} $W_x$ to the particles
$x\in \Z^d$
which determine how likely a particle may play the role of a hub\index{hub} in the resulting network.

\begin{figure}[htb!]
\begin{center}
\begin{minipage}[t]{0.45\textwidth}
\begin{center}
\includegraphics[width=.9\textwidth]{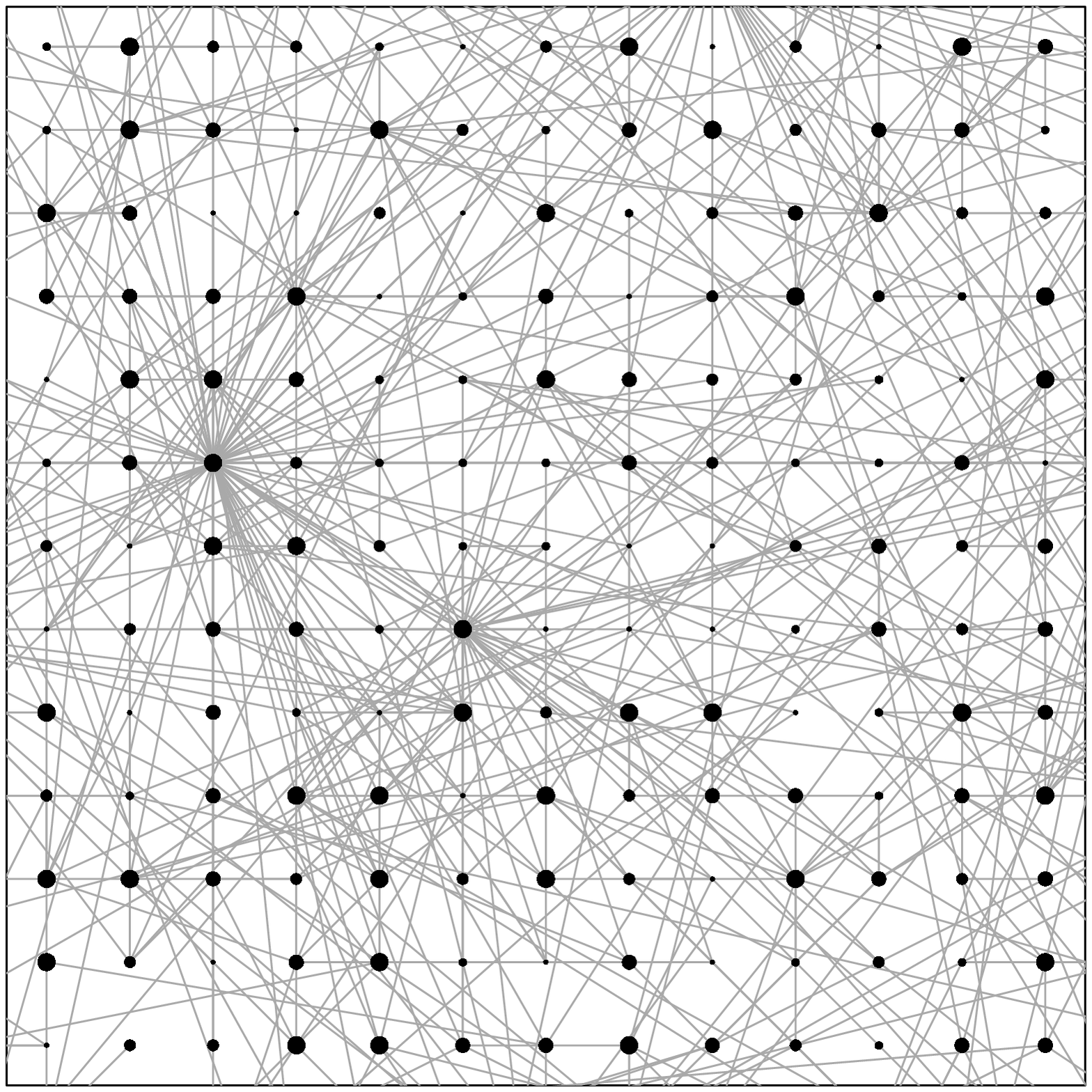}
\end{center}
\end{minipage}
\begin{minipage}[t]{0.45\textwidth}
\begin{center}
\includegraphics[width=.9\textwidth]{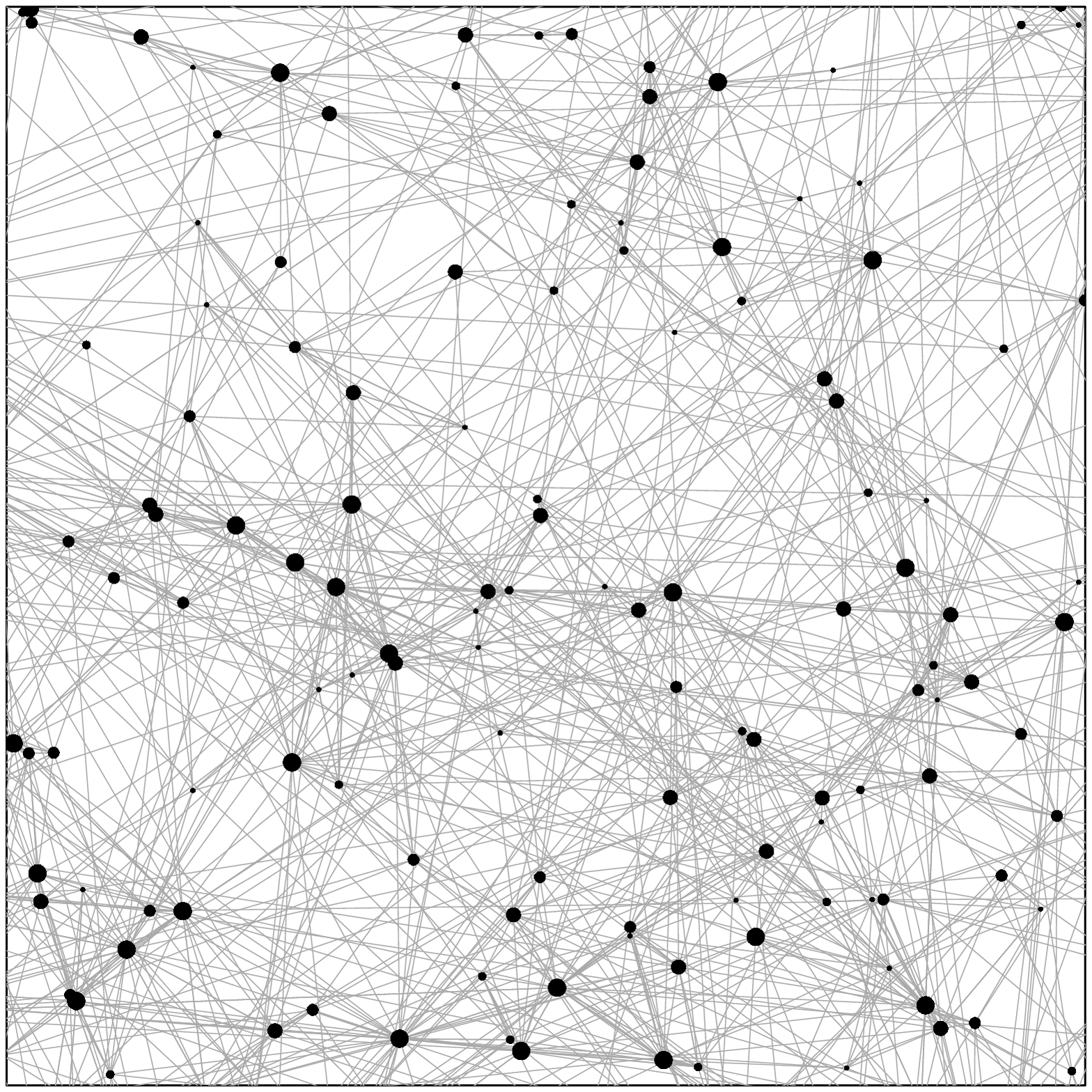}
\end{center}
\end{minipage}
\end{center}
\caption{{\bf lhs}: heterogeneous long-range percolation; 
{\bf rhs}: continuum space long-range percolation; the size of the 
particles illustrates the different weights $W_x\ge 1$} 
\label{HetLRP_and_RCM}
\end{figure}

Consider again the percolation model on the lattice $\Z^d$. Assume that $(W_x)_{x\in \Z^d}$
are i.i.d.~Pareto distributed with threshold parameter 1 and 
tail parameter $\beta>0$, i.e.~for $w\ge 1$
\begin{equation}\label{definition of Pareto}
\p\left[ W_x \le w \right] = 1 - w^{-\beta}.
\end{equation}
Choose $\alpha>0$ and $\lambda>0$ fixed. Conditionally given
$(W_x)_{x\in \Z^d}$, we consider the edge
probabilities\index{edge probability} for $x,y\in \Z^d$ given by
\begin{equation}\label{heterogeneous long-range percolation}
p_{x,y} =
1-\exp(-\lambda W_xW_y\|x-y\|^{-\alpha}).
\end{equation}
Between any pair $x,y\in \Z^d$ we attach an edge, independently of all other
edges, as follows
\begin{equation*}
\eta_{x,y}=\eta_{y,x} =
\left\{ 
\begin{array}{ll}
1 \qquad & \text{ with probability $p_{x,y}$,}\\
0 \qquad & \text{ with probability $1-p_{x,y}$.}
\end{array}
\right.
\end{equation*}
We denote the resulting probability 
measure on the edge configurations by $\p_{\lambda, \alpha, \beta}$. 
In contrast to (\ref{long-range percolation}) we have additional weights $W_x$ and $W_y$
in (\ref{heterogeneous long-range percolation}). The bigger these weights the more likely
is an edge between $x$ and $y$. Thus, particles $x\in \Z^d$ with a big weight $W_x$ will have
many adjacent particles $y$ (i.e.~particles $y\in \Z^d$ with $\eta_{x,y}=1$). 
Such particles $x$ will play the role of hubs in the network system.  
Figure \ref{HetLRP_and_RCM} (lhs) shows part of a realised edge configuration. 

The first interesting result is that this model provides a heavy-tailed degree
distribution, see Theorems 2.1 and 2.2 in Deijfen et al.~\cite{Remco}.
Denote again by ${\cal D}(0)$ the number of particles of $\Z^d$ that are adjacent to 0, then
we have the following result.
\begin{theorem} \label{infinite variance case}
Fix $d\ge 1$. We have the following two cases for the
degree distribution:
\begin{itemize}
\item for $\min\{\alpha, \beta \alpha\} \le d$, a.s., ${\cal D}(0)=\infty$;
\item for $\min\{\alpha, \beta \alpha\} > d$ set $\tau=\beta\alpha/d$, then
\begin{equation*}
\p_{\lambda, \alpha,\beta}\left[{\cal D}(0) > k \right]
=k^{-\tau} \ell(k), 
\end{equation*}
for some function $\ell(\cdot)$ that is slowly varying at infinity.
\end{itemize}
\end{theorem}
We observe that the heavy-tailedness of the weights $W_x$ induces heavy-tailedness 
in the degree distribution which is similar to choice (\ref{NSW degree})
in the NSW random graph model of Sect.~\ref{Section NSW graph}.
For $\alpha>d$ there are three different regimes: 
(i) $\beta\alpha \le d$ implies infinite degree, a.s.; (ii) for
$d<\beta\alpha <2d$ the degree distribution has finite mean but infinite variance
because $1<\tau<2$;
(iii) for
$\beta\alpha>2d$ the degree distribution has finite variance because $\tau>2$.
We will see that the distinction of the latter two cases has also implications on the
behaviour of the percolation properties and the graph distances similar to the considerations
in NSW random 
graphs. Note that from a practical point of view the interesting regime is (ii).

We again consider the connected component of a given particle $x\in \Z^d$
denoted by ${\cal C}(x)$ and we define the  percolation probability 
(for given $\alpha$ and $\beta$)
\begin{equation*}
\theta(\lambda)= \p_{\lambda, \alpha, \beta} \left[ |{\cal C}(0)|=\infty \right].
\end{equation*}
The {\it critical percolation value}\index{percolation!critical value} $\lambda_c$
is then defined by
\begin{equation*}
\lambda_c = \inf \left\{ \lambda>0:~\theta(\lambda)>0\right\}.
\end{equation*}
We have the following result, see Theorem 3.1 in Deijfen et al.~\cite{Remco}.
\begin{theorem} \label{theo 31}
Fix $d\ge 1$. Assume
$\min \{ \alpha, \beta \alpha \} > d$.
\begin{itemize}
\item[(a)] If $d\ge 2$, then $\lambda_c<\infty$.
\item[(b)]  If $d=1$ and $\alpha \in (1,2]$, then $\lambda_c<\infty$.
\item[(c)]  If $d=1$ and $\min \{ \alpha, \beta \alpha \} > 2$,
 then $\lambda_c=\infty$.
\end{itemize}
\end{theorem}
\noindent
This result is in line with Theorem \ref{theorem infinite cluster}.
Since $W_x\ge 1$, a.s., an edge configuration from edge probabilities $p_{x,y}$ defined in
(\ref{heterogeneous long-range percolation}) stochastically dominates 
an edge configuration with edge
probabilities $1- \exp(-\lambda \|x-y\|^{-\alpha})$. The latter is similar to the
homogeneous long-range percolation model on $\Z^d$ and
the results of the above theorem directly follow from Theorem \ref{theorem infinite cluster}.
For part (c) of the theorem we also refer to 
Theorem 3.1 of Deijfen et al.~\cite{Remco}. The next theorem follows
from Theorems 4.2 and 4.4 of Deijfen et al.~\cite{Remco}.

\begin{theorem} \label{theo 32}
Fix $d\ge 1$. Assume
$\min \{ \alpha, \beta \alpha \} > d$.
\begin{itemize}
\item[(a)] If $\beta \alpha <2d$, then $\lambda_c = 0$.
\item[(b)] If $\beta \alpha >2d$, then $\lambda_c > 0$.
\end{itemize}
\end{theorem}

\begin{figure}[htb!]
\begin{center}
\includegraphics[width=.65\textwidth]{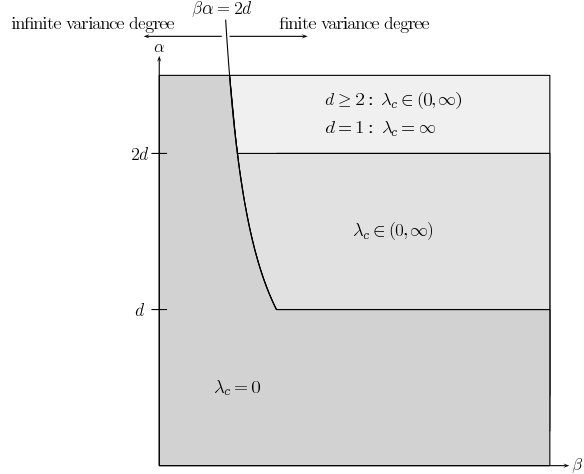}
\end{center}
\caption{phase transition picture for $d\ge1$}
\label{Picture: Phase}
\end{figure}

\noindent
Theorems \ref{theo 31} and \ref{theo 32} give the phase
transition\index{phase transition} pictures for $d\ge1$, see Fig.~\ref{Picture: Phase} for an illustration. 
They differ for $d=1$ and $d\ge 2$ in that
the former has a region where $\lambda_c=\infty$ and the latter
does not, see also the distinction in Theorem \ref{theorem infinite cluster}.
The most interesting case from a practical point of view is the infinite variance
case, $1<\tau<2$ and 
$d<\beta\alpha<2d$, respectively, which provides percolation for any $\lambda>0$.
It follows from Gandolfi et al.~\cite{GKN} that there is only   
one infinite connected component $\mathcal{C}_\infty$ whenever 
$\lambda>\lambda_c$, a.s. 
A difficult question to answer is what happens at criticality for $\lambda_c>0$. There
is the following partial result, see Theorem 3 in Deprez et al.~\cite{HW}: 
for $\alpha \in (d,2d)$ and $\beta \alpha > 2d$,
there does not exist an infinite connected component at criticality $\lambda_c>0$.
The case $\min\{\alpha,\beta\alpha\}>2d$ is still open.

\begin{figure}[htb!]
\begin{center}
\includegraphics[width=.9\textwidth]{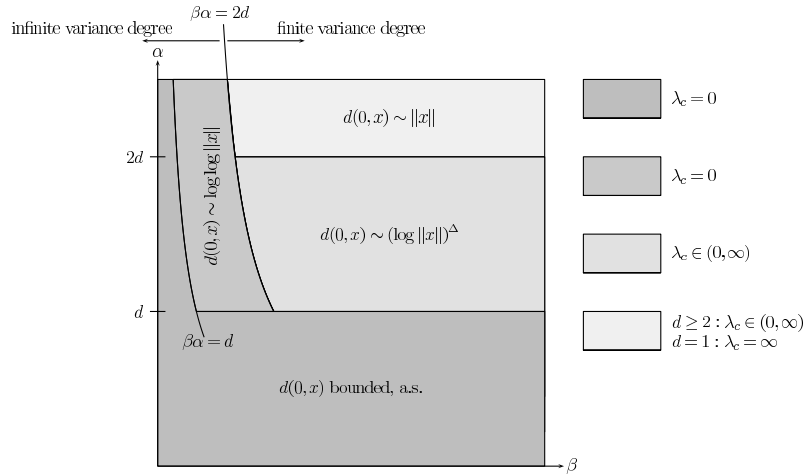}
\end{center}
\caption{chemical distances according to 
Theorem \ref{graph distance} and corresponding conjectures}
\label{Picture: Distances}
\end{figure}

Next we consider the graph distance $d(x,y)$, see also (\ref{definition of graph distance d}).
We have the following result, see Deijfen et al.~\cite{Remco}
and Theorem 8 in Deprez et al.~\cite{HW}. 
\begin{theorem}\label{graph distance}
Assume $\min \{ \alpha, \beta \alpha \} > d$.
\begin{itemize}
\item[(a)] (infinite variance of degree distribution $1<\tau<2$). Assume $d<\beta\alpha<2d$. For
any $\lambda>\lambda_c=0$ there exists $\eta_1> 0$ such that
for all $\epsilon>0$
\begin{equation*}
\lim_{\|x\|\to \infty}
\p_{\lambda, \alpha, \beta} \left[ \left. \eta_1\le \frac{d(0,x)}{\log \log \|x\|}\le 
 \frac{2}{|\log(\beta\alpha/d-1)|}+\epsilon
 \right|0,x\in {\cal C}_\infty\right]=1.
\end{equation*}
\item[(b1)] (finite variance of degree distribution $\tau>2$ case 1). Assume that $\beta\alpha >2d$
and $\alpha \in (d,2d)$. For any $\lambda>\lambda_c$ 
there exists $\eta_2\ge 1 $ such that 
for all $\epsilon>0$
\begin{equation*}
\lim_{\|x\|\to \infty}
\p_{\lambda, \alpha, \beta} \left[ \left. 
\eta_2-\epsilon
\le \frac{\log d(0,x)}{\log\log \|x\|}\le 
\Delta+\epsilon
\right|0,x\in {\cal C}_\infty\right]=1,
\end{equation*}
where $\Delta$ was defined in Theorem \ref{theorem homogeneous 1}.
\item[(b2)] (finite variance of degree distribution $\tau>2$ case 2). Assume 
 $\min\{\alpha,\beta\alpha\} >2d$. There exists $\eta_3>0$ such that
\begin{equation*}
\lim_{\|x\|\to \infty}
\p_{\lambda, \alpha, \beta} \left[\eta_3<  \frac{ d(0,x)}{ \|x\|} \right]=1.
\end{equation*}
\end{itemize}
\end{theorem}
Compare Theorem \ref{graph distance} (heterogeneous case) to
Theorems \ref{theorem homogeneous 1} and \ref{theorem homogeneous 2} (homogeneous case).
We observe that in the finite variance cases (b1)--(b2), i.e.~for $\tau=\beta\alpha/d>2$,
we obtain the same behaviour
for heterogeneous and homogeneous long-range percolation models. 
The infinite variance case (a) of the
degree distribution, i.e.~$1<\tau <2$ and $d<\beta\alpha<2d$, 
respectively, is new.
 This infinite variance case
provides a much slower decay of the graph distance, that is $d(0,x)$ is
of order $\log \log \|x\|$ as $\|x\|\to \infty$. This is a pronounced version
of the small-world effect, and this behaviour is similar to the NSW random graph model. 
Recall that empirical studies often suggest a tail
parameter $\tau$ between 1 and 2 which corresponds to the 
infinite variance regime of the degree distribution. 
In Fig.~\ref{Picture: Distances} we illustrate Theorem \ref{graph distance}
and we complete the picture 
about the chemical distances with the corresponding conjectures. 

We conclude that this model fulfils all three stylised facts of small-world effect,
the clustering property (which is induced by the Euclidean distance in the
probability weights (\ref{heterogeneous long-range percolation})) and the heavy-tailedness
of the degree distribution.

\section{Continuum space long-range percolation model\index{percolation!continuum space long-range}}
\label{homogeneous Poisson model}
The model of last section is restricted to the lattice $\Z^d$. A straightforward modification
is to replace the lattice $\Z^d$ by a homogeneous Poisson point process $X$ in $\R^d$. In comparison to
the lattice model, some of the proofs simplify because we can apply classical integration
in $\R^d$, other proofs become more complicated because one needs to make sure that the
realisation of the Poisson point process is sufficiently regular in space. As in Deprez-W\"uthrich
\cite{DW} we consider a homogeneous marked Poisson point process in $\R^d$, where
\begin{itemize}
\item $X$ denotes the spatially homogeneous Poisson point 
process\index{Poisson point process} in $\R^d$ with constant
intensity $\nu>0$. The individual particles of $X$ are denoted by $x\in X \subset \R^d$;
\item $W_x$, $x\in X$, are i.i.d.~marks\index{marks} having a Pareto distribution with
threshold parameter 1 and tail parameter $\beta>0$, see (\ref{definition of Pareto}).
\end{itemize}
Choose $\alpha>0$ and $\lambda>0$ fixed. Conditionally given $X$ and 
$(W_x)_{x\in X}$, we consider the edge
probabilities\index{edge probability} for $x,y\in X$ given by
\begin{equation}\label{Poisson long-range percolation}
p_{x,y} =
1-\exp(-\lambda W_xW_y\|x-y\|^{-\alpha}).
\end{equation}
Between any pair $x,y\in X$ we attach an edge, independently of all other
edges, as follows
\begin{equation*}
\eta_{x,y}=\eta_{y,x} =
\left\{ 
\begin{array}{ll}
1 \qquad & \text{ with probability $p_{x,y}$,}\\
0 \qquad & \text{ with probability $1-p_{x,y}$.}
\end{array}
\right.
\end{equation*}
We denote the resulting probability 
measure on the edge configuration by $\p_{\nu,\lambda, \alpha, \beta}$. 
Figure \ref{HetLRP_and_RCM} (rhs) shows part of a realised configuration. 
We have the following result for the degree distribution, see
Proposition 3.2 and Theorem 3.3 in Deprez-W\"uthrich \cite{DW}.
\begin{theorem} \label{infinite variance case Poisson}
Fix $d\ge 1$. We have the following two cases for the
degree distribution:
\begin{itemize}
\item for $\min\{\alpha, \beta \alpha\} \le d$, a.s., $\p_0[{\cal D}(0)=\infty|W_0]=1$;
\item for $\min\{\alpha, \beta \alpha\} > d$ set $\tau=\beta\alpha/d$, then
\begin{equation*}
\p_{0}\left[{\cal D}(0) > k \right]
=k^{-\tau} \ell(k),
\end{equation*}
for some function $\ell(\cdot)$ that is slowly varying at infinity.
\end{itemize}
\end{theorem}

\noindent
{\it Remarks.}
\begin{itemize}
\item Note that the previous statement needs some care because we need to 
make sure that there is a particle at the origin. This is not straightforward in the
Poisson case and $\p_0$ can be understood as the conditional distribution, conditioned
on having a particle at the origin. The formally precise construction is known as the
Palm distribution\index{Palm distribution}, 
which considers distributions shifted by the particles in the Poisson
cloud $X$.
\item In analogy to the homogeneous long-range percolation model in $\Z^d$ we
could also consider continuum space homogeneous long-range percolation in $\R^d$. This is achieved
by setting $W_x=W_y=1$, a.s., in (\ref{Poisson long-range percolation}). 
In this case the proof of the statement equivalent to (\ref{homogeneous lattice})
becomes rather easy. We briefly give the details in the next
lemma, see also proof of Lemma 3.1 in
Deprez-W\"uthrich \cite{DW}.
\end{itemize}

\begin{lemma} \label{lemma Poisson degrees}
Choose $W_x=W_y=1$, a.s., in (\ref{Poisson long-range percolation}).
For $\alpha\le d$ we have, a.s., ${\cal D}(0)=\infty$; for $\alpha>d$ the 
degree ${\cal D}(0)$ has a Poisson distribution.
\end{lemma}

\noindent
\begin{Proof}[of Lemma \ref{lemma Poisson degrees}
and (\ref{homogeneous lattice}) in continuum space]~
Let $X$ be a Poisson cloud with $0\in X$ and denote by 
$X(A)$ the number of particles in $X\cap A$ for $A\subset\R^d$. 
Every particle $x\in X\setminus\{0\}$ is
now independently from the others
removed from the Poisson cloud with probability $1-p_{0,x}$. 
The resulting process $\widetilde{X}$ is a thinned Poisson cloud 
having intensity function 
$x\mapsto\nu p_{0,x}=\nu(1-\exp(-\lambda\|x\|^{-\alpha}))\sim\nu\lambda\|x\|^{-\alpha}$
as $\|x\|\to \infty$. 
Since $\mathcal{D}(0)=\tilde X(\R^d\setminus\{0\})$ it follows that 
$\mathcal{D}(0)$ is infinite, a.s., if $\alpha\le d$ and that $\mathcal{D}(0)$ has 
a Poisson distribution otherwise. 
To see this let 
%
$\mu$ denote the Lebesgue measure in $\R^d$ and choose
a finite Borel set $A \subset \R^d$ containing the origin.
Since $A$ contains the origin, we have
$X(A) \ge \widetilde{X}(A)\ge 1$. This motivates
for $k\in \mathbb{N}_0$ to study
\begin{eqnarray*}
\p_0\left[\widetilde{X}(A)=k+1 \right]
&=& \sum_{i\ge k}
\p_0\left[\left.\widetilde{X}(A)=k+1\right|X(A)=i+1 \right]
\p_0\left[{X}(A)=i+1 \right].
\end{eqnarray*}
Since $A$ contains the origin, the case $i=0$ is trivial, 
i.e.~$\p_0[\widetilde{X}(A)=1 |X(A)=1 ]=1$. There remains $i\ge 1$.
Conditionally on $\{{X}(A)=i+1\}$, the $i$ particles (excluding the
origin) are independent and uniformly distributed in $A$. The
conditional moment generating function for $r\in \R$ is then
given by
\begin{eqnarray*}
&&\hspace{-1cm}
\EE_0 \left[\left. \exp \left\{ r\left(\widetilde{X}(A)-1\right)\right\}\right|X(A)=i+1 \right]
\\&&=~
\frac{1}{\mu (A)^i} \int_{A\times \cdots \times A}
\EE_0 \left[ \exp \left\{ r \sum_{l=1}^i \eta_{0,x_l}
\right\} \right]
dx_1\cdots dx_i
\\
&&=~
\frac{1}{\mu (A)^i} \int_{A\times \cdots \times A}
\prod_{l=1}^i
\EE_0 \left[ \exp \left\{ r \eta_{0,x_l}
\right\} \right]
dx_1\cdots dx_i
\\&&=~
\left(\frac{1}{\mu (A)} \int_{A}
\EE_0 \left[ \exp \left\{ r \eta_{0,x}
\right\} \right]
dx\right)^i.
\end{eqnarray*}
We calculate the integral for $W_0=W_{x}=1$, a.s., in (\ref{Poisson long-range percolation})
\begin{eqnarray*}
\frac{1}{\mu (A)} \int_{A}
\EE_0 \left[ \exp \left\{ r \eta_{0,x}
\right\} \right] dx
&=&
\frac{1}{\mu (A)} \int_{A}
 e^r p_{0,x}+ (1- p_{0,x})~
dx
\\&=&
e^rp(A)
+ \left(1- p(A)\right),
\end{eqnarray*}
with $p(A)=\mu (A)^{-1} \int_{A}
p_{0,x} dx \in (0,1)$.
Thus, conditionally on $\{X(A)=i+1\}$,
$\widetilde{X}(A)-1$ has a binomial distribution with parameters $i$ and $p(A)$.
This implies that
\begin{eqnarray*}
\p_0\left[\widetilde{X}(A)=k+1 \right]
&=&\sum_{i\ge k}
\binom{i}{k}~ p(A)^k \left(1-p(A)\right)^{i-k}~
\p_0\left[{X}(A)=i+1 \right]
\\&=&\sum_{i\ge k}
\frac{p(A)^k}{k!} \frac{\left(1-p(A)\right)^{i-k}}{(i-k)!}~
\exp \{-\nu \mu(A)\}(\nu \mu(A))^i
\\&=&\frac{(\nu \mu(A)p(A))^k}{k!} 
\sum_{i\ge k}
\frac{\left(\nu \mu(A)(1-p(A))\right)^{i-k}}{(i-k)!}
\exp \{-\nu \mu(A)\}
\\&=&
\exp \left\{-\nu \mu(A)p(A)\right\}
\frac{(\nu \mu(A)p(A))^k}{k!}
\\&=&
\exp \left\{-\nu\int_{A}
p_{0,x} dx\right\}
\frac{\left( \nu\int_{A}
p_{0,x} dx\right)^k}{k!}.
\end{eqnarray*}
This implies that $\widetilde{X}$ is a non-homogeneous
Poisson point process with 
intensity function 
\begin{equation*}
x\mapsto \nu p_{0,x}
=\nu \left(1-\exp(-\lambda \|x\|^{-\alpha} )\right)
\sim \nu \lambda \|x\|^{-\alpha}, \qquad \text{ as $\|x\|\to \infty$.}
\end{equation*}
But this immediately implies that the degree distribution
${\cal D}(0)=\widetilde{X}(\R^d\setminus\{0\})$ is infinite, a.s., if $\alpha \le d$, 
and that it has a Poisson distribution otherwise. This finishes the proof.
\qed
\end{Proof}

~

\noindent
We now switch back to the heterogeneous long-range percolation
model (\ref{Poisson long-range percolation}).
We consider the connected component $\mathcal{C}(0)$ of a particle in the origin
under the Palm distribution $\p_0$.
We define the  percolation probability\index{percolation!probability} 
\begin{equation*}
\theta(\lambda)= \p_{0} \left[ |{\cal C}(0)|=\infty \right]. 
\end{equation*}
The {\it critical percolation value}\index{percolation!critical value}  $\lambda_c$
is then defined by
\begin{equation*}
\lambda_c = \inf \left\{ \lambda>0:~\theta(\lambda)>0\right\}.
\end{equation*}
We have the following results, see Theorem 3.4 in Deprez-W\"uthrich \cite{DW}.
\begin{theorem} \label{theo 31 DW}
Fix $d\ge 1$. Assume
$\min \{ \alpha, \beta \alpha \} > d$.
\begin{itemize}
\item[(a)] If $d\ge 2$, then $\lambda_c<\infty$.
\item[(b)]  If $d=1$ and $\alpha \in (1,2]$, then $\lambda_c<\infty$.
\item[(c)]  If $d=1$ and $\min \{ \alpha, \beta \alpha \} > 2$,
 then $\lambda_c=\infty$.
\end{itemize}
\end{theorem}

\begin{theorem} \label{theo 32 DW}
Fix $d\ge 1$. Assume
$\min \{ \alpha, \beta \alpha \} > d$.
\begin{itemize}
\item[(a)] If $\beta \alpha <2d$, then $\lambda_c = 0$.
\item[(b)] If $\beta \alpha >2d$, then $\lambda_c > 0$.
\end{itemize}
\end{theorem}
These are the continuum space analogues to Theorems \ref{theo 31} and \ref{theo 32},
for an illustration see also Fig.~\ref{Picture: Phase}.
The work on the graph distances  in the continuum space long-range percolation model
is still work in progress, but we expect similar results  to the ones in Theorem \ref{graph distance},
see also Fig.~\ref{Picture: Distances}.
However, proofs in the continuum space model are more sophisticated due to the
randomness of the positions of the particles.

The advantage of the latter continuum space model (with homogeneous marked Poisson point process)
is that it can be extended to non-homogeneous Poisson point processes. For instance, if certain
areas are more densely populated than others we can achieve such a non-homogeneous space
model by modifying the
constant intensity $\nu$ to a space-dependent density function $\nu(\cdot):\R^d \to \R_+$.

\section{Renormalisation techniques}
\label{Renormalisation techniques}
\begin{figure}
\begin{center}
\includegraphics[width=.6\textwidth]{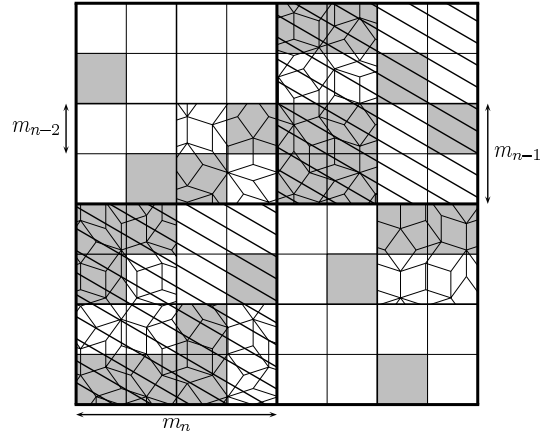}
\end{center}
\vspace{-.5cm}
\caption{
Example of the renormalisation technique. Define inductively 
the box lengths $m_n=2m_{n-1}$, $n\in\mathbb{N}$, for some initial $m_0\in\mathbb{N}$. 
Call translates of $[0,m_{0}-1]^d$ to be  $0$-stage boxes and assume 
that goodness of such boxes is defined. 
For $n\in\mathbb{N}$ we call the translates of $[0,m_{n}-1]^d$ $n$-stage boxes 
and we inductively say that an $n$-stage box is good if it contains at least two 
good $(n-1)$-stage boxes. 
Assume that in the illustration good $(n-2)$-stage boxes, $n\ge2$, 
are coloured in grey. The good $(n-1)$-stage boxes 
are then the Penrose-patterned boxes of side-length $m_{n-1}$ 
and the good $n$-stage boxes 
are striped. The illustrated $(n+1)$-stage box is good, because it 
contains two good $n$-stage boxes 
}
\label{Figure: renormalisation}
\end{figure}

In this section we present a crucial technique that is used in many of the proofs of the
previous statements. These proofs are often based on renormalisation techniques. 
That is, one collects particles
in boxes. These boxes are defined to be either {\it good} (having a certain property) or {\it bad}
(not possessing this property). These boxes are then again merged to bigger good or bad boxes. These 
scalings\index{scaling} and renormalisations are done over several generations of box  sizes, 
see Fig.~\ref{Figure: renormalisation} for an illustration.  
The purpose
of these rescalings is that one arrives at a certain generation of box sizes 
that possesses certain characteristics
to which classical site-bond percolation results apply. 
We exemplify this with a particular example. 

\subsection{Site-bond percolation\index{percolation!site-bond} }
Though we will not directly use site-bond percolation, we start with the
description of this model because it is often useful.
Site-bond percolation in $\Z^d$ is a modification of homogeneous long-range percolation
introduced in Sect.~\ref{homogeneous long-range percolation}.
Choose a fixed dimension $d\ge 1$ and consider the square lattice $\Z^d$. 
Assume that every site $x\in \Z^d$ is occupied independently with probability $r^\ast\in [0,1]$ and 
every bond between $x$ and $y$ in $\Z^d$ is occupied independently with 
probability
\begin{equation}\label{site-bond percolation}
p^\ast_{x,y} =
1-\exp(-\lambda^\ast \|x-y\|^{-\alpha}),
\end{equation}
for given parameters $\lambda^\ast>0$ and $\alpha>0$.
The connected component ${\cal C}^\ast(x)$ of a given site $x\in \Z^d$ is then defined to be the
set of all occupied sites $y\in \Z^d$ such that $x$ and $y$ are connected by a path only running through
occupied sites and occupied bonds (if $x$ is not occupied then ${\cal C}^\ast(x)$ is the empty set).
We can interpret this as follows: we place particles at sites $x\in \Z^d$ at random with
probability $r^\ast$. This defines a (random) subset of $\Z^d$ and then we consider long-range percolation
on this random subset, i.e.~this corresponds to a thinning of homogeneous
long-range percolation in $\Z^d$. 
We can then study the percolation properties of this site-bond percolation model, some
results are presented in Lemma 3.6 of Biskup \cite{Biskup} and
in the proof of Theorem 2.5 of Berger \cite{Berger1}.
The aim in many proofs in percolation theory
is to define different generations of box sizes using renormalisations, see Fig.~\ref{Figure: renormalisation}.
We perform these renormalisations until we arrive at a generation
of box sizes for which good boxes occur sufficiently often. 
If this is the case and if all the necessary dependence assumptions
are fulfilled we can apply classical site-bond percolation results. 

In order
to simplify our outline we use a modified version of the homogeneous long-range 
percolation\index{percolation!modified homogeneous long-range}  
model (\ref{long-range percolation}) of Sect.~\ref{homogeneous long-range percolation}.
We set $p=1-\exp(-\lambda)$ and obtain the following model.

\begin{model}[modified homogeneous long-range percolation] 
\label{long-range percolation model}
Fix $d\ge 2$.
Choose $\alpha>0$ and  $\lambda>0$ fixed and define the
edge probabilities\index{edge probability}  for $x,y \in \Z^d$ by
\begin{equation*}
p_{x,y} =
1-\exp(-\lambda \|x-y\|^{-\alpha}).
\end{equation*}
Then edges between all 
pairs of particles
$x,y \in \Z^d$ are attached independently with edge probability $p_{x,y}$
and the probability measure of the resulting edge configurations $\eta=(\eta_{x,y})_{x,y\in \Z^d}$
is denoted by $\p_{\lambda, \alpha}$.
\end{model}
Note that this model is a special case of site-bond percolation with $r^\ast=1$ and
$\lambda^\ast=\lambda$ in (\ref{site-bond percolation}).

\subsection{Largest semi-clusters\index{semi-cluster}}

In order to demonstrate the renormalisation technique we repeat the
proof of  Lemma 2.3 of Berger \cite{Berger1}
in the modified homogeneous long-range percolation Model \ref{long-range percolation model}, 
see Theorem \ref{connected component} below. This proof is rather sophisticated because it needs
a careful treatment of dependence and we revisit the second version of the 
proof of Lemma 2.3 provided in Berger \cite{Berger11}.

Fix $\alpha\in(d,2d)$ and choose $\lambda>0$ so large 
that there exists a unique infinite connected  
component, a.s., having density $\kappa>0$
(which exists due to Theorem \ref{theorem infinite cluster}).
Choose $M \ge 1$ and $K\ge 0$ integer valued.
For $v\in \Z^d$ we define box\index{box} $B_v$ and its $K$-enlargement\index{$K$-enlargement} $B^{(K)}_v$ by
\begin{equation*}
B_v = Mv + [0,M-1]^d \qquad \text{ and } \qquad
B^{(K)}_v = Mv + [-K,M+K-1]^d.
\end{equation*}
For every box $B_v$ we define a $\ell$-\textit{semi-cluster} to be 
a set of at least $\ell$ sites in $B_v$ which are connected within $B_v^{(K)}$. 
For any $\varepsilon>0$ there exists $M'\ge 1$ such that for all $M\ge M'$
and some $K\ge 0$ we have
\begin{eqnarray}\label{good probability}
\p_{\lambda, \alpha} [\text{at least $M^d\kappa/2$ sites of $B_v$ belong to the infinite connected }
~
\\\text{component and these sites are connected within $B^{(K)}_v$}] &\ge& 1- \varepsilon/2.
\nonumber
\end{eqnarray}
Existence of $M'\ge 1$ follows from the ergodic theorem and existence of $K$ from the fact
that the infinite connected component is unique, a.s., and therefore all sites in $B_v$ belonging
to the infinite connected component need to be connected within a certain $K$-enlargement of $B_v$.
Formula (\ref{good probability}) says that we have a $(M^d\kappa/2)$-semi-cluster in $B_v$
with at least probability $1-\varepsilon/2$. We first show uniqueness of large
semi-clusters.

\begin{lemma} \label{uniqueness of semi-clusters}
Choose $\xi \in (\alpha/d, 2)$ and $\gamma \in (0,1)$ with $18\gamma>16+\xi$.
There exist $\varphi=\varphi(\xi,\gamma)>0$ 
and $M'=M'(\xi,\gamma)\ge 1$ such that for all $M\ge M'$ and all $K\ge 0$
we have
\begin{equation*}
\p_{\lambda, \alpha}
\left[\text{there is at most one 
$M^{d\gamma}$-semi-cluster in $B_v$}\right] > 1-M^{-d\varphi},
\end{equation*}
where by ``at most one'' we mean that there is no second $M^{d\gamma}$-semi-cluster in $B_v$
which is not connected to the first one within $B_v^{(K)}$.
\end{lemma}

\noindent
\begin{Proof}[of Lemma \ref{uniqueness of semi-clusters}]
The proof uses the notion of inhomogeneous random 
graphs\index{random graph!inhomogeneous} as defined in Aldous \cite{Aldous}. 
An inhomogeneous random graph $H(N,\xi)$ with size $N$ and 
parameter $\xi$ is a set of particles $\{1,\ldots,k\}$ and corresponding 
masses $s_1,\ldots,s_k$ such that $N=\sum_{i=1}^ks_i$; and any  $i\neq j$ are connected 
independently with probability $1-\exp\left(-s_is_jN^{-\xi}\right)$. 
From Lemma 2.5 of Berger \cite{Berger11} we know that for any $1<\xi<2$ and  
$0<\gamma<1$ with $18\gamma>16+\xi$, there exist $\varphi=\varphi(\xi,\gamma)>0$ 
and $N'=N'(\xi,\gamma)\ge 1$ such that for all $N\ge N'$ and every inhomogeneous 
random graph with size $N$ and parameter $\xi$ we have 
\begin{eqnarray}\label{Equation: Lemma 2.5}
&&\hspace{-.5cm}\nonumber
\p\left[
\text{$H(N,\xi)$ contains more than one connected component $C$ with 
$\sum_{i\in C}s_i \ge N^\gamma$}
\right]
\\&&<~ N^{-\varphi}.
\end{eqnarray}
We now show uniqueness of $M^{d\gamma}$-semi-clusters in $B_v$. 
Choose $\xi \in (\alpha/d, 2)$ and $\gamma \in (0,1)$ such that  $18\gamma>16+\xi$.
For any $x,y\in B_v$ we have $\|x-y\| \le \sqrt{d}M$. Choose $M$ so large that
$\lambda(\sqrt{d}M)^{-\alpha}>M^{-d\xi}$ and choose $K\ge 0$ arbitrarily. 
Particles $x,y\in B_v$ are then attached with
probability $p_{x,y}$ uniformly bounded by
\begin{eqnarray*}
p_{x,y} &=&
1-\exp(-\lambda \|x-y\|^{-\alpha})
\\&\ge&
1- \exp\left(-\lambda\left(\sqrt{d}M\right)^{-\alpha}\right)
> 
1- \exp\left(-M^{-d\xi}\right) = \nu>0,
\end{eqnarray*}
where the last equality defines $\nu$. This allows to decouple the sampling
of edges $\eta=(\eta_{x,y})_{x,y\in \Z^d}$ in $B_v$.
For every $x,y\in B_v$, define $p'_{x,y}\in (0,1)$ by
\begin{equation*}
p_{x,y}= p'_{x,y} + \nu - \nu p'_{x,y}.
\end{equation*} 
We now sample $\eta=(\eta_{x,y})_{x,y\in\Z^d}$ in two steps. 
We first sample $\eta'$ according to Model \ref{long-range percolation model} but
with edge probabilities $p'_{x,y}$ if $x,y \in B_v$ and with 
edge probabilities $p_{x,y}$ otherwise. 
Secondly, we sample $\eta''$ as an independent configuration on $B_v$ 
where there is an edge between $x$ and $y$ with edge probability $\nu$ for  
$x,y \in B_v$. 
By definition of $p'_{x,y}$ we get that $\eta'\vee\eta'' \stackrel{\rm (d)}{=}\eta$. 
Let $S_{1},S_2 \subset B_v$ be two disjoint maximal sets of sites in 
$B_v$ that are $\eta'$-connected within $B_v^{(K)}$, i.e.~$S_1$ and $S_2$ are two disjoint
maximal semi-clusters in $B_v$ for given edge configuration $\eta'$.  
Note that by maximality 
\begin{eqnarray*}
&&\hspace{-1cm}
\p\left[\left.
\text{there is an $\eta$-edge between $S_1$ and $S_2$}\right|\eta'
\right]
\\&&=
\p\left[\left.
\text{there is an $\eta''$-edge between $S_1$ and $S_2$}\right|\eta'
\right]
\\&&=1-(1-\nu)^{|S_{1}||S_{2}|}~=~
1- \exp\left(-|S_{1}||S_{2}|M^{-d\xi}\right). 
\end{eqnarray*}
If we denote by $S_1, \ldots, S_k$ all disjoint maximal semi-clusters in $B_v$
for given edge configuration $\eta'$ then we see that these maximal semi-clusters
form an inhomogeneous random graph of size $\sum_{i=1}^k |S_i| =M^d$ and
parameter $\xi$. Therefore, there exist $\varphi>0$ and $M'\ge 1$ such that for all
$M\ge M'$ and all $K\ge 0$
we have from (\ref{Equation: Lemma 2.5})
\begin{eqnarray*}
&&\hspace{-.1cm}
\p\Bigg[
\text{$H(M^d,\xi)$ contains more than one connected component $C$}
~
\\&&\hspace{6cm}\text{with 
$\sum_{i\in C}|S_i| \ge M^{d\gamma}$}\Bigg|\eta'
\Bigg]
<M^{-d\varphi}.
\end{eqnarray*}
Note that this bound is uniform in $\eta'$ and $K\ge 0$. Therefore, the probability of having at least two
$M^{d\gamma}$-semi-clusters in $B_v$ which are not connected within $B_v^{(K)}$ is bounded
by $M^{-d\varphi}$. \qed
\end{Proof}

~

\noindent
We can now combine (\ref{good probability}) and Lemma \ref{uniqueness of semi-clusters}.
Choose $\varepsilon>0$. For all $M$ sufficiently large and $K\ge 0$ such that
(\ref{good probability}) holds we have
\begin{equation}\label{exaclty one large semi-cluster}
\p_{\lambda, \alpha}
\left[\text{there is exactly one 
$(M^{d}\kappa/2)$-semi-cluster in $B_v$}\right] \ge 1-\varepsilon,
\end{equation}
where by ``exactly one'' we mean that there is no other 
$(M^{d}\kappa/2)$-semi-cluster in $B_v$ which is not connected to the first
one within $B_v^{(K)}$. This follows because of $\gamma<1$, which implies
that $M^{d\gamma}\le M^d \kappa/2$ for all $M$ sufficiently large, and
because $M^{-d\varphi}<\varepsilon/2$ for all $M$ sufficiently large.

\subsection{Renormalisation}
Choose $\varepsilon>0$ fixed, 
and $M> 1$ and $K\ge 0$ such that (\ref{exaclty one large semi-cluster}) holds. For $v\in \Z^d$ we
say that box $B_v$ is {\it good} if there is exactly one $(M^d\kappa/2)$-semi-cluster
in $B_v$ (where exactly one is meant in the sense of above). 
Therefore, on good boxes there are at least 
$M^d\kappa/2$ sites in $B_v$ that are connected within $B^{(K)}_v$ 
and we have
\begin{equation}\label{good probability 2}
\p_{\lambda, \alpha} [ \text{$B_v$ is good} ] \ge 1- \varepsilon.
\end{equation}
Note that the goodness properties of $B_{v_1}$ and $B_{v_2}$ for $v_1\neq v_2\in \Z^d$ are not necessarily
independent because their $K$-enlargements
$B^{(K)}_{v_1}$ and $B^{(K)}_{v_2}$ may overlap.

Now, we define renormalisation over different generations $n\in \mathbb{N}_0$; terminology  
 {\it $n$-stage}\index{box!$n$-stage} is referred to the $n$-th generation.
Choose an integer valued sequence $a_n>1$, $n\in \mathbb{N}_0$, 
with $a_0=M$ and define the box lengths
$(M_n)_{n\in \mathbb{N}_0}$ as follows: set $M_0=a_0=M$ and for $n\in\mathbb{N}$
\begin{equation*}
M_n= a_nM_{n-1}=M_0\prod_{i=1}^n a_i=\prod_{i=0}^n a_i.
\end{equation*}
Define the $n$-stage boxes, $n\in \mathbb{N}_0$, by
\begin{equation*}
B_{n,v} = M_nv + [0,M_n-1]^d \qquad \text{ with $v\in \Z^d$.}
\end{equation*}
Note that $n$-stage boxes $B_{n,v}$ have volume
$M^d_n= a_n^dM^d_{n-1}=\prod_{i=0}^n a^d_i$ and every $n$-stage box $B_{n,v}$
contains $a_n^d$ of $(n-1)$-stage boxes $B_{n-1,x}\subset B_{n,v}$, and
$(M_n/a_0)^d= \prod_{i=1}^n a^d_i$
of $0$-stage boxes $B_{x}=B_{0,x}\subset B_{n,v}$, see also Fig.~\ref{Figure: renormalisation}.

~

\noindent
{\bf Renormalisation.\index{renormalisation}}
We define goodness of $n$-stage boxes $B_{n,v}$ recursively for a given sequence
$\kappa_n \in (0,1)$, $n\in \mathbb{N}_0$, of densities where we initialise
$\kappa_0 = \kappa/2$.

~

\noindent
{\it (i) Initialisation $n=0$.} We say that $0$-stage box $B_{0,v}$, $v\in \Z^d$,
is good\index{box!good} if it contains exactly one $(\kappa_0a_0^d)$-semi-cluster.
Due to our choices of $M> 1$ and $K\ge 0$ we see that the goodness
of $0$-stage box $B_{0,v}$ occurs with at least probability $1-\varepsilon$,
see (\ref{good probability 2}). 

~

\noindent
{\it (ii) Iteration $n-1 \to n$.} Choose $n\in \mathbb{N}$ and assume that goodness
of $(n-1)$-stage boxes $B_{n-1,v}$, $v\in \Z^d$, has been defined. For $v\in \Z^d$ 
we say that $n$-stage box $B_{n,v}$ is good if the event 
$A_{n,v}=
A^{(a)}_{n,v}\cap A^{(b)}_{n,v}$ occurs, where 
\begin{itemize}
 \item[(a)] $A^{(a)}_{n,v}=\{\text{at least $\kappa_{n}a_{n}^d$ of the $(n-1)$-stage boxes
$B_{n-1,x}\subset B_{n,v}$ are good}\}$; and
\item[(b)] $A^{(b)}_{n,v}=\{\text{all 
$(\prod_{i=0}^{n-1}\kappa_ia_i^d)$-semi-clusters of all good 
$(n-1)$-stage boxes }$\\
$\text{\hspace{6cm} in $B_{n,v}$ are connected within $B^{(K)}_{n,v}$}\}$. 
\end{itemize}
 \hfill  \qed

~

\noindent
Observe that on event $A_{n,v}$ the $n$-stage box $B_{n,v}$ 
contains at least $\prod_{i=0}^{n}\kappa_ia_i^d$ sites that 
are connected within the $K$-enlargement $B^{(K)}_{n,v}$ of $B_{n,v}$.  
We set density $u_n=\prod_{i=0}^n\kappa_i$ which gives  
\begin{equation*}
\prod_{i=0}^n\kappa_i a_i^d=M_n^d\prod_{i=0}^n\kappa_i=M_n^d u_n. 
\end{equation*}
Therefore, good $n$-stage boxes contain $(M_n^d u_n)$-semi-clusters.
Our next aim is to calculate the probability $p_n$ of having a good $n$-stage box.
The case $n=0$ follows from (\ref{good probability 2}), i.e.~for any
$\varepsilon>0$ and any $M$ sufficiently large there exists $K\ge 0$ 
such that
\begin{equation*}
p_0~=~\p_{\lambda, \alpha} [ \text{$B_{0,v}$ is good} ]
~=~\p_{\lambda, \alpha} [ \text{$B_v$ is good} ] ~\ge~ 1-\varepsilon.
\end{equation*}




\begin{theorem}\label{connected component}
Assume $\alpha \in (d,2d)$. Choose $\lambda>0$ so large that we have a unique
infinite connected component, a.s., having density $\kappa>0$.
For every $\varepsilon' \in (0,1)$ there exists $N_0\ge 1$
such that for all $N\ge N_0$
\begin{equation*}
\p_{\lambda, \alpha}\left[|C_N| \ge  N^{\alpha/2}\right]\ge 1-\varepsilon',
\end{equation*}
where $C_N$ is the largest connected component\index{connected component!largest} in $[0,N-1]^d$.
\end{theorem}
Note that for density $\kappa>0$ of the infinite connected component we expect
roughly $\kappa N^d$ sites in box $[0,N-1]^d$ belonging to the infinite connected
component. The above lemma however says
that at least $N^{\alpha/2}$ sites in $[0,N-1]^d$ are connected {\it within that
box}. That is, here we do not need any $K$-enlargements as in (\ref{good probability})
and, therefore, this event is independent for different
disjoint boxes $vN+[0,N-1]^d$ and we may apply classical 
site-bond percolation results.

~

\noindent
\begin{Proof}[of Theorem \ref{connected component}]
Choose $\alpha\in(d, 2d)$ and $\varepsilon' \in (0,1)$ fixed. As in Lemma 2.3
of Berger \cite{Berger11} we now make a choice of parameters and sequences which will provide
the statement of Theorem \ref{connected component}.
Choose $\xi \in (\alpha/d, 2)$ and $\gamma \in (0,1)$ such that  $18\gamma>16+\xi$.
Choose $\delta>\vartheta>1$ with $2\vartheta<\delta(2d-\alpha)$
and $d\delta-\vartheta>d\gamma \delta$.
Note that this is possible because it requires
that $\delta \min\{1,(2d-\alpha)/2,d(1-\gamma)\}>\vartheta>1$.
 Define for $n\in \mathbb{N}$ 
\begin{equation}\label{sequences}
\kappa_n=(n+1)^{-\vartheta} \qquad \text{ and } \qquad a_n=(n+1)^\delta.
\end{equation}
For simplicity, we assume that $\delta$ is an integer which implies that also
$a_n>1$ is integer valued, and $\kappa_n \in (0,1)$ will play the role of densities introduced
above.
Observe that for $\vartheta >1$ we have for all $n\ge 1$
\begin{equation}\label{c_1}
\prod_{l=1}^n (1+3\kappa_l)~ \le~ \lim_{n\to \infty}\prod_{l=1}^n (1+3\kappa_l)~=~c_1 \in (1,\infty).
\end{equation}
Choose $\varepsilon \in (0,\varepsilon'/c_1)\subset(0,1)$ fixed. 

There still remains the choice of $a_0=M\ge 1$ and $\kappa_0\in (0,1)$. We set $\kappa_0=\kappa/2$.
Note that choices (\ref{sequences}) imply 
\begin{equation*}
(2M_{n-1}^d)^{\gamma}=2^\gamma M^{d\gamma}
(n!)^{d\gamma\delta} \qquad \text{ and } \qquad
M_{n-1}^d u_{n-1}=\frac{\kappa}{2}M^d 
(n!)^{d\delta-\vartheta}.
\end{equation*}
Therefore,
\begin{equation}\label{rhs uniformly bounded in n}
\frac{M_{n-1}^d u_{n-1}}{(2M_{n-1}^d)^{\gamma}}=\frac{\kappa}{2^{1+\gamma}}M^{d(1-\gamma)} 
(n!)^{d\delta-\vartheta-d\gamma\delta}.
\end{equation}
Because of $d\delta-\vartheta>d\gamma\delta$ the right-hand side of
(\ref{rhs uniformly bounded in n}) is uniformly bounded from below in $n\ge 1$
and for $M$ sufficiently large the right-hand side of (\ref{rhs uniformly bounded in n})
is strictly bigger than 1 for all $n\ge 1$. Therefore, there exists $m_1\ge 1$
such that for all $M\ge m_1$ and all $n\ge 1$ we have
\begin{equation}\label{crucial inequality 1}
(2M_{n-1}^d)^{\gamma}<M_{n-1}^d u_{n-1}.
\end{equation}
Next we are going to bound for $n\in \mathbb{N}_0$ the probabilities
\begin{equation*}
p_n~=~\p_{\lambda, \alpha} [ \text{$B_{n,v}$ is good} ]
~=~\p_{\lambda, \alpha} [ A_{n,v} ].
\end{equation*}
We have for $n\ge 1$
\begin{eqnarray}
1-p_n&=&\p_{\lambda, \alpha} \left[A^{\rm c}_{n,v}\right]\nonumber
\\&=& \p_{\lambda, \alpha} \left[(A^{(a)}_{n,v}\cap A^{(b)}_{n,v})^{\rm c}\right]
\label{bound one}
\le  \p_{\lambda, \alpha} \left[(A^{(a)}_{n,v})^{\rm c} \right]
+\p_{\lambda, \alpha}\left[ (A^{(b)}_{n,v})^{\rm c}\right].
\end{eqnarray}
For the first term in (\ref{bound one}) we have,
using Markov's inequality and translation invariance,
\begin{eqnarray*}
\p_{\lambda, \alpha} \left[
(A^{(a)}_{n,v})^{\rm c}\right]
&=&
\p_{\lambda, \alpha} \left[\sum_{B_{n-1,x}\subset B_{n,v}}
1_{A_{n-1,x}} < \kappa_{n}a_{n}^d \right]
\\&=&
\p_{\lambda, \alpha} \left[\sum_{B_{n-1,x}\subset B_{n,v}}
1_{A_{n-1,x}^{\rm c}} > (1-\kappa_{n})a_{n}^d \right]
\\&\le&
\frac{1}{(1-\kappa_{n})a_{n}^d}~
\sum_{B_{n-1,x}\subset B_{n,v}}
\p_{\lambda, \alpha} \left[A_{n-1,x}^{\rm c}\right]
\\&=&
\frac{1}{1-\kappa_{n}}~
\p_{\lambda, \alpha} \left[A_{n-1,v}^{\rm c}\right]~=~\frac{1-p_{n-1}}{1-\kappa_{n}}.
\end{eqnarray*}
The second term in (\ref{bound one}) is more involved due
to possible dependence in the $K$-enlargements. 
Choose $\varphi=\varphi(\xi,\gamma)>0$ and $M'(\xi,\gamma)\ge
1$ as in Lemma \ref{uniqueness of semi-clusters}.
On event $(A^{(b)}_{n,v})^{\rm c}$ there exist at least two $(M_{n-1}^d u_{n-1})$-semi-clusters 
in good $(n-1)$-stage boxes $B_{n-1,v_1}$ and $B_{n-1,v_2}$ in
$B_{n,v}$ that are not connected within the $K$-enlargement $B_{n,v}^{(K)}$.
Define $B=B_{n-1,v_1} \cup B_{n-1,v_2}$. Note that $B$
has volume $2M_{n-1}^d$ and that any $x,y \in B_{n,v}$ have maximal distance $\sqrt{d}M_n$.
We analyse the following ratio
\begin{equation*}
\frac{\lambda(\sqrt{d}M_n)^{-\alpha}}{(2M^d_{n-1})^{-\xi}}
=\lambda d^{-\alpha/2}2^{\xi}~M^{d \xi-\alpha}
(n!)^{(d \xi-\alpha)\delta} ~(n+1)^{-\alpha \delta}.
\end{equation*}
Note that $d\xi>\alpha$. This implies that
the right-hand side of the previous equality is uniformly
bounded from below in $n\ge 1$. Therefore, there exists $m_2\ge m_1$ such that for all
$M\ge m_2$ and all $n\ge 1$ inequality (\ref{crucial inequality 1}) holds and
\begin{equation} \label{rhs uniformly bounded in n 2}
\lambda(\sqrt{d}M_n)^{-\alpha}>(2M^d_{n-1})^{-\xi}.
\end{equation}
This choice implies that for any $x,y \in B$ we have
\begin{eqnarray*}
p_{x,y} &=&
1-\exp(-\lambda \|x-y\|^{-\alpha})
\\&\ge&
1- \exp\left(-\lambda\left(\sqrt{d}M_n\right)^{-\alpha}\right)
> 
1- \exp\left(-(2M^d_{n-1})^{-\xi}\right) = \nu_n>0,
\end{eqnarray*}
where the last equality defines $\nu_n$.
We now proceed as in Lemma \ref{uniqueness of semi-clusters}.
Decouple the sampling
of edges $\eta=(\eta_{x,y})_{x,y\in \Z^d}$ in $B$.
For every $x,y\in B$, define $p'_{x,y}\in (0,1)$ by
\begin{equation*}
p_{x,y}= p'_{x,y} + \nu_n - \nu_n p'_{x,y}.
\end{equation*} 
We again sample $\eta=(\eta_{x,y})_{x,y\in\Z^d}$ in two steps. 
We first sample $\eta'$ according to Model \ref{long-range percolation model} but
with edge probabilities $p'_{x,y}$ if $x,y \in B$ and with 
edge probabilities $p_{x,y}$ otherwise. 
Secondly, we sample $\eta''$ as an independent configuration on $B$ 
where there is an edge between $x$ and $y$ with edge probability $\nu_n$ for  
$x,y \in B$. 
By definition of $p'_{x,y}$ we get that $\eta'\vee\eta'' \stackrel{\rm (d)}{=}\eta$. 
Let $S_{1},S_2 \subset B$ be two disjoint maximal sets of sites in 
$B$ that are $\eta'$-connected within $B_{n,v}^{(K)}$, i.e.~$S_1$ and $S_2$ are two disjoint
maximal semi-clusters in $B$ for given edge configuration $\eta'$.  
Note that by maximality 
\begin{eqnarray*}
&&\hspace{-1cm}
\p\left[\left.
\text{there is an $\eta$-edge between $S_1$ and $S_2$}\right|\eta'
\right]
\\&&=
\p\left[\left.
\text{there is an $\eta''$-edge between $S_1$ and $S_2$}\right|\eta'
\right]
\\&&=~1-(1-\nu_n)^{|S_{1}||S_{2}|}~=~
1- \exp\left(-|S_{1}||S_{2}|(2M^d_{n-1})^{-\xi}\right). 
\end{eqnarray*}
If we denote by $S_1, \ldots, S_k$ all disjoint maximal semi-clusters in $B$
for given edge configuration $\eta'$ then we see that these maximal semi-clusters
form an inhomogeneous random graph of size $\sum_{i=1}^k |S_i| =2M_{n-1}^d \ge 2 M^d$ and
parameter $\xi$. Therefore, for  choices $\varphi=\varphi(\xi,\gamma)>0$ 
and $m_3 \ge \max\{m_2, M'(\xi,\gamma)\}$ (where $\varphi(\xi,\gamma)$ and $M'(\xi,\gamma)$ were given
by Lemma \ref{uniqueness of semi-clusters}) 
we have that for all $M\ge m_3$ and
all $n\ge 1$ inequality (\ref{crucial inequality 1})  holds,
and for all $K\ge 0$
we have from (\ref{Equation: Lemma 2.5})
\begin{eqnarray*}
&&\hspace{-.2cm}
\p\Bigg[
\text{$H(2M_{n-1}^d,\xi)$ contains more than one connected component $C$ }
~
\\&&\hspace{5.2cm}\text{with 
$\sum_{i\in C}|S_i| \ge (2M_{n-1}^d)^{\gamma}$}\Bigg|\eta'
\Bigg]
< (2M_{n-1}^d)^{-\varphi}.
\end{eqnarray*}
Note that this bound is uniform in $\eta'$ and $K\ge 0$ and holds for all $n\ge 1$. 
Therefore, the probability of having at least two
$(2M_{n-1}^d)^{\gamma}$-semi-clusters in $B$ which are not connected within 
$B_{n,v}^{(K)}$ is bounded by $(2M_{n-1}^d)^{-\varphi}$. 
Next we use that for all $m\ge m_3$ inequality (\ref{crucial inequality 1})
holds. Therefore, we get for  all $M\ge m_3$,
all $n\ge 1$ and all $K\ge 0$
\begin{equation*}
 \p_{\lambda, \alpha} \left[
\text{there are at least two $(M_{n-1}^d u_{n-1})$-semi-clusters in $B$}\right]
< (2M_{n-1}^d)^{-\varphi}.
\end{equation*}
Note that $B_{n,v}$ contains $a_n^d$ disjoint $(n-1)$-stage boxes, therefore we get
for all $M\ge m_3$, all $n\ge 1$ and all $K\ge 0$
\begin{equation*}
\p_{\lambda,\alpha}\left[ (A^{(b)}_{n,v})^{\rm c}\right]
~\le~ 
\binom{a_n^d}{2}\left(2M_{n-1}^{d}\right)^{-\varphi}
~\le~
a_n^{2d}M_{n-1}^{-d\varphi}. 
\end{equation*}
This implies for  all $M\ge m_3$,
all $n\ge 1$ and all $K\ge 0$
\begin{eqnarray*}
1-p_n&=&\p_{\lambda, \alpha} [ \text{$B_{n,v}$ is not good} ]
\\&=&\p_{\lambda, \alpha} [ A_{n,0}^c ]\le \frac{1-p_{n-1}}{1-\kappa_{n}}+a_n^{2d}M_{n-1}^{-d\varphi}
\le (1-p_{n-1})(1+2\kappa_n)+a_n^{2d}M_{n-1}^{-d\varphi}.
\end{eqnarray*}
Consider
\begin{equation*}
\frac{a_n^{2d}M_{n-1}^{-d\varphi}}{\varepsilon \kappa_n}=
\frac{(n+1)^{2d\delta}M^{-d\varphi}(n!)^{-d\delta\varphi}}{\varepsilon 
(n+1)^{-\vartheta}}
=\varepsilon^{-1}M^{-d\varphi}
(n+1)^{2d\delta+\vartheta}(n!)^{-d\delta\varphi}.
\end{equation*}
Note that this is uniformly bounded from above in $n$. Therefore, there exists
$m_4\ge m_3$ such that for  all $M\ge m_4$,
all $n\ge 1$ and all $K\ge 0$
\begin{equation*}
1-p_n \le (1-p_{n-1})(1+2\kappa_n)+\varepsilon \kappa_n
\le 
(1+3\kappa_n)\max\{1-p_{n-1},\varepsilon\}.
\end{equation*}
Applying induction we obtain for all $M\ge m_4$,
all $n\ge 1$ and all $K\ge 0$
\begin{equation*}
1-p_n 
\le 
\max\{\varepsilon, 1-p_0\} \prod_{i=1}^n(1+3\kappa_i).
\end{equation*}
Choose $m_5 \ge m_4$ such that for all $M\ge m_5$ there exists $K=K(M)\ge 0$ such
that (\ref{exaclty one large semi-cluster})  and (\ref{good probability 2}) hold.
These choices imply that $p_0\ge 1-\varepsilon$. 
Therefore, for all $M\ge m_5$, $K(M)$ such that (\ref{good probability 2}) holds, and 
all $n\ge 1$
\begin{equation*}
1-p_n 
\le \varepsilon \prod_{i=1}^n(1+3\kappa_i)\le 
\varepsilon c_1
< \varepsilon',
\end{equation*}
where $c_1\in (1,\infty)$ was defined in (\ref{c_1}).
Thus, for all $M\ge m_5$,
$K(M)$ such that (\ref{good probability 2}) holds, and for all $n\ge 0$
\begin{eqnarray}\label{lower bound provides}
&&\hspace{-1cm}
\p_{\lambda,\alpha}\left[\text{there are at least $M_n^du_n$ sites in box $B_{n,0}$ connected within $B_{n,0}^{(K)}$}\right]
\\&&\ge~
\p_{\lambda,\alpha} [A_{n,0}] \ge 1- \varepsilon',\nonumber
\end{eqnarray}
note that $B_{n,0}=[0,M_n-1]^d$.
Note that the explicit choices (\ref{sequences}) provide
\begin{equation*}
M_n =M ((n+1)!)^\delta \qquad \text{ and } \qquad
u_n = \kappa_0 ((n+1)!)^{-\vartheta}.
\end{equation*}
The edge length of $B_{n,0}^{(K)}$ is given by
$M_n + 2 K = M_n + 2 K(M) = 
M ((n+1)!)^\delta + 2 K(M)$. Therefore, for all $M\ge 1$
there exists $n_0\ge 1$ such that for all $n\ge n_0$
\begin{equation*}
M_n + 2 K  \le N=N(M,n)=2 M ((n+1)!)^\delta,
\end{equation*}
where the last identity is the definition of $N=N(M,n)$.
For all $n \ge n_0$
the number of connected vertices in $B_{n,0}^{(K)}$ under (\ref{lower bound provides})
is at least
\begin{equation*}
M^d_n u_n = \kappa_0 M^d((n+1)!)^{-\vartheta+d\delta}
=2^{\vartheta/\delta-d}\kappa_0M^{\vartheta/\delta}N^{d-\vartheta/\delta}.
\end{equation*}
Note that the choices of $\delta$ and $\vartheta$ are such that
$d-\vartheta/\delta>\alpha/2>0$. Therefore, there exists $n_1\ge n_0$ such that
for all $n\ge n_1$
we have
\begin{equation*}
M^d_n u_n 
\ge N^{\alpha/2}.
\end{equation*}
This implies for all $n\ge n_1$, see (\ref{lower bound provides}),
\begin{eqnarray*}
&&\hspace{-.5cm}
\p_{\lambda,\alpha}\left[|C_N| \ge  N^{\alpha/2}\right]
\\&&\ge
\p_{\lambda,\alpha}\left[\text{there are at least $N^{\alpha/2}$ sites in box $B_{n,0}$ connected within $B_{n,0}^{(K)}$}\right]
\\&&\ge
\p_{\lambda,\alpha}\left[\text{there are at least $M_n^du_n$ sites in box $B_{n,0}$ connected within $B_{n,0}^{(K)}$}\right]
\\&&\ge~1- \varepsilon',
\end{eqnarray*}
where $N=N(M,n)\ge N(M,n_1)=2M((n_1+1)!)^\delta$.
This  proves the claim on the grid
$N(M,n_1),N(M,n_1+1),\ldots, $ with
$N(M,n+1)=N(M,n)(n+2)^\delta$ for $n\ge n_1$.
For $n' \in [N(M,n), N(M,n+1))$ we have 
on the set $\{|C_{N(M,n)}| \ge \rho_0 N(M,n)^{d-\vartheta/\delta}\}$
with $\rho_0=2^{\vartheta/\delta-d}\kappa_0M^{\vartheta/\delta}$
\begin{eqnarray*}
|C_{n'}| &\ge& |C_{N(M,n)}| ~\ge~ \rho_0 N(M,n)^{d-\vartheta/\delta}
\\&=&
\rho_0 N(M,n)^{d-\vartheta/\delta-\alpha/2}
\left(\frac{N(M,n)}{n'}\right)^{\alpha/2}
(n')^{\alpha/2}
\\&\ge&\rho_0 N(M,n)^{d-\vartheta/\delta-\alpha/2}
\left(\frac{N(M,n)}{N(M,n+1)}\right)^{\alpha/2}
(n')^{\alpha/2}
\\&=&\rho_0 ~\frac{N(M,n)^{d-\vartheta/\delta-\alpha/2}}
{\left(n+2\right)^{\delta\alpha/2}}~
(n')^{\alpha/2}~\ge~ (n')^{\alpha/2},
\end{eqnarray*}
for all $n$ sufficiently large.
This finishes the proof of Theorem \ref{connected component}.\qed
\end{Proof} 

~

\noindent
{\bf Conclusion.} Theorem \ref{connected component} defines good boxes $[0,N-1]^d$ on a new scale, i.e.~these
are boxes that contain sufficiently large connected components $C_N$. The latter occurs with probability 
$1-\varepsilon'\ge r^\ast$, for small $\varepsilon'$. If we can prove that such large connected components
in disjoint boxes are connected
by an occupied edge with probability bounded below by (\ref{site-bond percolation}),
then we are in the set-up of a site-bond percolation model. This is exactly what is used in Theorem 3.2
of Biskup \cite{Biskup} in order to prove that (i) large connected components are percolating, a.s.;
and (ii) $|C_N|$ is even of order $\rho N^{d}$ for an appropriate positive constant $\rho>0$, which
improves Theorem \ref{connected component}.

\bibliography{References}
\end{document}